\documentclass[11pt]{article}
\usepackage{epsfig,graphicx,ifpdf}

\def\be{\begin{equation}}
\def\ee{\end{equation}}
\def\bea{\begin{eqnarray}}
\def\eea{\end{eqnarray}}
\def\bes{\begin{eqnarray*}}
\def\ees{\end{eqnarray*}}

\def\nn{\nonumber}
\def\lb{\label}
\def\bs{\setminus}

\def\vs{{\varsigma}}

\def\R{{\bf R}}
\def\C{{\bf C}}
\def\Z{{\bf Z}}
\def\K{{\bf K}}
\def\N{{\bf N}}
\def\U{{\bf U}}

\def\Q{{\bf Q}}
\def\T{{\bf T}}
\def\CP{{\bf CP}}

\def\aa{{\alpha}}
\def\bb{{\beta}}
\def\ga{{\gamma}}

\def\ka{{\kappa}}
\def\th{{\theta}}
\def\om{{\omega}}
\def\Om{{\Omega}}
\def\ep{{\epsilon}}
\def\lm{{\lambda}}
\def\Lm{{\Lambda}}

\def\sg{{\sigma}}
\def\dm{{\diamond}}

\def\vf{{\varphi}}

\def\S{{\cal S}}
\def\K{{\cal K}}
\def\P{{\cal P}}

\def\E{{\cal E}}

\def\tdR{\tilde{R}}

\def\rank{{\rm rank}}
\def\Sp{{\rm Sp}}
\def\mod{{\rm mod}}
\def\ol{\overline}

\def\hb{\vrule height0.18cm width0.14cm $\,$}

\title{The index growth and multiplicity of closed geodesics}

\author{Huagui Duan$^{1}$\thanks{Partially supported by NNSF Grant 10801079,
RFDP Grant 200800551002, LPMC of MOE of China and Nankai University.
E-mail: duanhg@nankai.edu.cn} \quad and \quad Yiming
Long$^{2}$\thanks{Partially supported by the 973 Program of MOST,
Yangzi River Professorship, NNSF, MCME, RFDP, LPMC of MOE of China,
and Nankai University. E-mail: longym@nankai.edu.cn}\\\\
$^{1}$ School of Mathematics\\
$^{2}$ Chern Institute of Mathematics and LPMC\\
Nankai University, Tianjin 300071\\ The People's Republic of China\\}
\date{}

\begin{document}

\maketitle

\begin{abstract}
{\it In the recent paper \cite{LoD1}, we classified closed geodesics on
Finsler manifolds into rational and irrational two families, and
gave a complete understanding on the index growth properties of iterates
of rational closed geodesics. This study yields that a rational closed
geodesic can not be the only closed geodesic on every irreversible or
reversible (including Riemannian) Finsler sphere, and that there exist
at least two distinct closed geodesics on every compact simply connected
irreversible or reversible (including Riemannian) Finsler $3$-dimensional
manifold. In this paper, we study the index growth properties of irrational
closed geodesics on Finsler manifolds. This study allows us to extend
results in \cite{LoD1} on rational and in \cite{DuL1}, \cite{Rad4} and \cite{Rad5}
on completely non-degenerate closed geodesics on spheres and $\CP^2$ to every
compact simply connected Finsler manifold. Then we prove the existence of at
least two distinct closed geodesics on every compact simply connected irreversible
or reversible (including Riemannian) Finsler $4$-dimensional manifold. }
\end{abstract}



\renewcommand{\theequation}{\thesection.\arabic{equation}}
\renewcommand{\thefigure}{\thesection.\arabic{figure}}

\setcounter{equation}{0}
\section{Introduction and main results}

It has been a long-standing problem in dynamical systems and differential
geometry whether every compact Riemannian manifold has infinitely many
distinct closed geodesics. D. Gromoll and W. Meyer \cite{GrM1} in 1969
proved the following result:

{\bf Theorem A.} (cf. \cite{GrM1}) {\it On a compact Riemannian manifold
there exist infinitely many closed geodesics, if the free loop space of
this manifold has an unbounded sequence of Betti numbers.}

Stimulated by this result, M. Vigu\'e-Poirrier and D. Sullivan \cite{ViS1}
in 1976 established following result:

{\bf Theorem B.} (cf. \cite{ViS1}) {\it The free loop space of a
compact simply connected Riemannian manifold $M$ has no unbounded sequence
of Betti numbers if and only if the rational cohomology algebra of $M$
possess only one generator. }

Both of the two theorems were generalized to corresponding Finsler manifolds by
H. Matthias in 1980 (cf. \cite{Mat1}). Therefore based on these two theorems,
the most interesting manifolds in this multiplicity problem are those
compact simply connected manifolds satisfying
\be  H^*(M;\Q)\cong T_{d,h+1}(x)=\Q[x]/(x^{h+1}=0)   \lb{1.1}\ee
with a generator $x$ of degree $d\ge 2$ and hight $h+1\ge 2$. The main
examples are the compact rank one symmetric spaces, i.e., spheres $S^d$
of dimension $d$ with $h=1$, complex projective spaces $\C P^h$ of dimension
$2h$ with $d=2$, quaternionic projective spaces $\mathbf{H}P^h$ of dimension
$4h$ with $d=4$, and the Cayley plane $\C aP^2$ of dimension $16$ with $d=8$
and $h=2$.

The studies of closed geodesics on such manifolds can be chased back
to J. Jacobi, J. Hadamard, H. Poincar\'e, G. D. Birkhoff, M. Morse,
L. Lyusternik and Schnirelmann and others. Specially G. D. Birkhoff
established the existence of at least one closed geodesic on every
Riemannian sphere $S^d$ with $d\ge 2$ (cf. \cite{Bir1}). Later L.
Lyusternik and A. Fet proved the existence of at least one closed
geodesic on every compact Riemannian manifold (cf. \cite{LyF1}).
An important breakthrough on this problem is due to V. Bangert \cite{Ban2}
and J. Franks \cite{Fra1} around 1990, who proved that there exist always
infinitely many closed geodesics on every Riemannian $2$-sphere (cf.
also \cite{Hin1} and \cite{Hin2}). But when the dimension of a compact simply
connected manifold is greater than $2$, we are not aware of any multiplicity
results on the existence of at least two closed geodesics without pinching or
bumpy conditions even on spheres (cf. \cite{Ano1}, \cite{Ban1}, \cite{Kli1},
\cite{BTZ1}, \cite{BTZ2} and \cite{DuL1}, \cite{Rad4}, \cite{Rad5}), except
the Theorem C below proved recently in \cite{LoD1}.

When one considers irreversible Finsler metrics, the problem of
counting closed geodesics becomes more delicate because of A. Katok's
famous example of 1973 which shows that there exist some irreversible Finsler
metrics on $S^d$ with only finitely many closed geodesics (cf. \cite{Kat1}
and \cite{Zil2}). In \cite{HWZ1} of 2003, H. Hofer, K. Wysocki and E. Zehnder
proved that there exist either two or infinitely many distinct prime closed
geodesics on a Finsler $(S^2,F)$ provided that all the iterates of all closed
geodesics are non-degenerate and the stable and unstable manifolds of all
hyperbolic closed geodesics intersect transversally. In \cite{BaL1} of 2005,
V. Bangert and Y. Long proved that on every irreversible Finsler $S^2$ there
always exist at least two distinct prime closed geodesics.

Note that in the recent \cite{LoD1}, we have proved the following

{\bf Theorem C.} {\it There exist always at least two distinct prime (geometrically
distinct) closed geodesics for every irreversible (or reversible, specially
Riemannian) Finsler metric on any $3$-dimensional compact simply connected manifold,
where the typical case is $S^3$.}

To further our study on the multiplicity of closed geodesics, we note that
in the famous book \cite{Mor1} of 1934, M. Morse studied closed
geodesics on ellipsoids. Specially he proved that for any given
integer $N>0$, every closed geodesic $c$ of a $d$-dimensional ellipsoid $E^d$
in $\R^{d+1}$ which is not an iterate of some main ellipse must have Morse
index $i(c)\ge N$, provided all the semi-axes of $E^d$ are less than $1$ and
sufficiently closed to $1$. Consequently the Betti numbers at dimensions
less than $N$ of the free loop space of such an $E^d$ can be generated by
iterates of the $d+1$ main ellipses on $E^d$ only. His this result suggests
that it is necessary to study asymptotic and growth properties of Morse
indices of iterates of prime closed geodesics on the manifold in order to
get multiplicity results.

In the recent paper \cite{LoD1}, we classified prime closed geodesics
on any compact Finsler manifold $M$ into two families: rational and
irrational. Here a prime closed geodesic is {\bf rational}, if its basic
normal form decomposition (cf. Section 3 below) introduced by Y. Long
in \cite{Lon1} and \cite{Lon2} contains no $2\times 2$ rotation matrix
$R(\th)=\left(\matrix{\cos\th & -\sin\th \cr
                      \sin\th & \cos\th \cr}\right)$ with $\th/\pi\in \R\bs\Q$,
and {\bf irrational} otherwise. A prime closed geodesic is {\bf completely
non-degenerate}, if all of its iterates $c^m$ are non-degenerate.

Recall that on a compact Finsler manifold $(M,F)$, a closed geodesic
$c:S^1=\R/\Z\to M$ is {\bf prime}, if it is not a multiple covering
(i.e., iteration) of any other closed geodesics. Here the $m$-th iteration
$c^m$ of $c$ is defined by $c^m(t)=c(mt)$ for $m\in\N$. The inverse curve
$c^{-1}$ of $c$ is defined by $c^{-1}(t)=c(1-t)$ for $t\in S^1$. Two
prime closed geodesics $c_1$ and $c_2$ on a Finsler manifold $(M,F)$
(or Riemannian manifold $(M,g)$) are {\bf distinct} (or {\bf
geometrically distinct}), if they do not differ by an $S^1$-action
(or $O(2)$-action).

In \cite{LoD1}, the index growth properties of rational closed geodesics
are completely understood. This result is used to prove that on every
(irreversible or reversible) Finsler sphere $S^d$, it is impossible that
there exists only one prime closed geodesic which is rational.

In the Section 3 of this paper, we study first the growth properties of
indices of iterates of irrational closed geodesics. We show that if the
initial index of a prime closed geodesic is not too small, then the Morse
indices $i(c^m)$ is monotone in $m\ge 1$. When this index monotonicity does
not hold, we prove that for a closed geodesic $c$, there exist infinitely
many positive integers $T$ such that the indices $\{i(c^m)\}_{m>T}$ and the
indices $\{i(c^m)\}_{m<T}$ are suitably separated by the sum of $i(c^T)$ and
some constant (see Theorem 3.21 below). We call this property the
{\bf quasi-monotonicity}.

As applications of these studies, in Section 4 we then generalize the
result in \cite{LoD1} on rational closed geodesics on spheres, and the
results in \cite{DuL1}, \cite{Rad4}, and \cite{Rad5} on completely
non-degenerate closed geodesics on spheres and $\C P^2$ to all compact
simply connected manifolds. That is:

\smallskip

{\bf Theorem 1.1.} {\it For every irreversible (or reversible, specially
Riemannian) Finsler metric $F$ on any compact simply connected manifold,
if there exists only one prime (geometrically distinct) closed geodesic,
it can be neither rational nor completely non-degenerate. }

Then using above results we study the $4$-dimensional case in the Sections
5 and 6 respectively, and prove the following theorems.

\smallskip

{\bf Theorem 1.2.} {\it For every irreversible Finsler metric $F$ on
any compact simply connected $4$-dimensional manifold, there always
exist at least two distinct prime closed geodesics.}

\smallskip

{\bf Theorem 1.3.} {\it For every reversible Finsler metric $F$ on
any compact simply connected $4$-dimensional manifold, there always
exist at least two geometrically distinct closed geodesics. In particular,
it holds for every such Riemannian manifold.}

\smallskip

For reader's conveniences, in Section 2 we briefly review some known results
on closed geodesics, and compute the precise sums of Betti numbers of the
$S^1$-invariant free loop space of compact simply connected manifolds
satisfying the condition (\ref{1.1}).

In this paper, we denote by $\N$, $\N_0$, $\Z$, $\Q$, $\R$, and $\C$ the
sets of positive integers, non-negative integers, integers, rational numbers,
real numbers, and complex numbers respectively. We define the functions
$[a]=\max\{k\in\Z\,|\,k\le a\}$, $\{a\}=a-[a]$,
$E(a)=\min\{k\in\Z\,|\,k\ge a\}$ and $\vf(a)=E(a)-[a]$. Denote by $\,^{\#}A$ the
number of elements in a finite set $A$. In this paper, we use only singular
homology modules with $\Q$-coefficients.

\setcounter{equation}{0}
\section{Critical point theory of closed geodesics}

\subsection{Critical modules for closed geodesics}

Let $M$ be a compact and simply connected manifold with a Finsler
metric $F$. Closed geodesics are critical points of the energy
functional $E(\gamma)=\frac{1}{2}\int_{S^1}F(\gamma(t),\dot{\gamma}(t))^2dt$
on the Hilbert manifold $\Lambda M$ of $H^1$-maps from $S^1$ to $M$. An
$S^1$-action is defined by $(s\cdot\gamma)(t)=\gamma(t+s)$ for all
$\gamma\in\Lm M$ and $s, t\in S^1$. The index form of the functional
$E$ is well defined along any closed geodesic $c$ on $M$, which we
denote by $E''(c)$. As usual, denote by $i(c)$ and $\nu(c)$ the
Morse index and nullity of $E$ at $c$. For a closed geodesic $c$,
denote by $c^m$ the $m$-fold iteration of $c$ and
$\Lm(c^m)=\{\ga\in\Lm M\,|\, E(\ga)<E(c^m)\}$. Recall that
respectively the {\it mean index} $\hat{i}(c)$ and the $S^1$-{\it
critical modules} of $c^m$ are defined by
\be \hat{i}(c)=\lim_{m\rightarrow\infty}\frac{i(c^m)}{m},
   \quad \ol{C}_*(E,c^m) = H_*\left((\Lm(c^m)\cup S^1\cdot c^m)/S^1,\Lm(c^m)/S^1\right).
        \lb{2.1}\ee

If $c$ has multiplicity $m$, then the subgroup $\Z_m=\{\frac{n}{m}:0\le n<m\}$
of $S^1$ acts on $\ol{C}_k(E,c)$. As on p.59 of \cite{Rad2}, for $m\ge1$, let
$H_*(X,A)^{\pm \Z_m}=\{[\xi]\in H_*(X,A):T_*[\xi]=\pm\xi\}$, where $T$ is a
generator of the $\Z_m$ action. On $S^1$-critical modules of $c^m$, the following
lemma holds:

\smallskip

{\bf Lemma 2.1.} (cf. Satz 6.11 of \cite{Rad2}) {\it Suppose $c$ is
a prime closed geodesic on a compact Finsler manifold $M$. Then
there exist two sets $U_{c^m}^-$ and $N_{c^m}$, the so-called local negative
disk and the local characteristic manifold at $c^m$ respectively,
such that $\nu(c^m)=\dim N_{c^m}$ and
\bea \overline{C}_q( E,c^m)
&\equiv& H_q\left((\Lm(c^m)\cup S^1\cdot c^m)/S^1, \Lm(c^m)/S^1\right)\nn\\
&=& \left(H_{i(c^m)}(U_{c^m}^-\cup\{c^m\},U_{c^m}^-)
    \otimes H_{q-i(c^m)}(N_{c^m}^-\cup\{c^m\},N_{c^m}^-)\right)^{+\Z_m}, \nn\eea

(i) When $\nu(c^m)=0$, there holds
$$ \overline{C}_q( E,c^m) = \left\{\matrix{
     \Q, &\quad {\it if}\;\; i(c^m)=i(c)\,(\mod 2)\;\;{\it and}\;\;
                   q=i(c^m),\;  \cr
     0, &\quad {\it otherwise}, \cr}\right.  $$

(ii) When $\nu(c^m)>0$, let $\ep(c^m)=(-1)^{i(c^m)-i(c)}$, then
there holds
$$ \overline{C}_q( E,c^m)=H_{q-i(c^m)}(N_{c^m}^-\cup\{c^m\},N_{c^m}^-)^{\ep(c^m)\Z_m}. $$}

Let
\be  k_j(c^m) \equiv \dim\, H_j( N_{c^m}^-\cup\{c^m\},N_{c^m}^-), \quad
     k_j^{\pm 1}(c^m) \equiv \dim\, H_j(N_{c^m}^-\cup\{c^m\},N_{c^m}^- )^{\pm\Z_m}.
        \lb{2.2}\ee
Then we have

\smallskip

{\bf Lemma 2.2.} (cf. \cite{Rad2}, \cite{BaL1}, \cite{DuL1}) {\it Let $c$
be a closed geodesic on a Finsler manifold $M$.

(i) There hold $0\le k_j^{\pm 1}(c^m) \le k_j(c^m)$ for $m\ge 1$ and
$j\in\Z$, $k_j(c^m)=0$ whenever $j\not\in [0,\nu(c^m)]$ and
$k_0(c^m)+k_{\nu(c^m)}(c^m)\le1$. If $k_0(c^m)+k_{\nu(c^m)}(c^m)=1$,
then $k_j(c^m)=0$ when $j\in(0,\nu(c^m))$.

(ii) For any $m\in\N$, there hold $k_0^{+1}(c^m) = k_0(c^m)$ and
$k_0^{-1}(c^m) = 0$. In particular, if $c^m$ is non-degenerate,
there hold $k_0^{+1}(c^m) = k_0(c^m)=1$, and $k_0^{-1}(c^m) =
k_j^{\pm 1}(c^m)=0$ for all $j\neq 0$.

(iii) Suppose for some integer $m=np\ge 2$ with $n$ and $p\in\N$ the
nullities satisfy $\nu(c^m)=\nu(c^n)$. Then there hold
$k_j(c^m)=k_j(c^n)$ and ${k}_j^{\pm 1}(c^m)={k}^{\pm 1}_j(c^n)$ for
any integer $j$.}

\subsection{Rademacher-type mean index identity for closed geodesics}

Let $(M,F)$ be a compact and simply connected Finsler manifold with
finitely many prime closed geodesics. It is well known that for every
prime closed geodesic $c$ on $(M,F)$, there holds either
$\hat{i}(c)>0$ and then $i(c^m)\to +\infty$ as $m\to +\infty$, or
$\hat{i}(c)=0$ and then $i(c^m)=0$ for all $m\in\N$. Denote those
prime closed geodesics on $(M,F)$ with positive mean indices by
$\{c_j\}_{1\le j\le k}$. In \cite{Rad1} and \cite{Rad2}, Rademacher
established a celebrated mean index identity relating all the $c_j$s
with the global homology of $M$ (cf. Section 7, specially Satz 7.9
of \cite{Rad2}) for compact simply connected Finsler manifolds.

For each $m\in\N$, let $\ep=\ep(c^m)=(-1)^{i(c^m)-i(c)}$ and
\bea  K(c^m)
&\equiv&(k_0^\ep(c^m), k_1^\ep(c^m), \ldots, k_{2\dim M-2}^\ep(c^m))\nn\\
&=&(k_0^{\ep (c^m)}(c^m),  k_1^{\ep(c^m)}(c^m), \ldots,
     k_{\nu (c^m)}^{\ep (c^m)}(c^m), 0, \ldots, 0). \lb{2.3}\eea

{\bf Lemma 2.3.} (cf. Lemmas 7.1 and 7.2 of \cite{Rad2}, cf. also \cite{LoD1})
{\it Let $c$ be a prime orientable closed geodesic on a compact Finsler manifold $(M,F)$.
Then there exist a minimal integer $N=N(c)\in\N$ such that $\nu(c^{m+N})=\nu(c^m)$,
$i(c^{m+N})-i(c^m)\in 2\Z$, and $K(c^{m+N})=K(c^m)$, $\forall\,m\in\N$.}

\smallskip

{\bf Lemma 2.4.} (Satz 7.9 of \cite{Rad2}, cf. also \cite{LoD1}) {\it Let
$(M,F)$ be a compact simply connected Finsler manifold with
$\,H^{\ast}(M,\Q)=T_{d,h+1}(x)$. Denote prime closed geodesics on $(M,F)$
with positive mean indices by $\{c_j\}_{1\le j\le k}$ for some $k\in\N$.
Then the following identity holds
\be \sum_{j=1}^k\frac{\hat{\chi}(c_j)}{\hat{i}(c_j)}=B(d,h)
=\left\{\matrix{
     -\frac{h(h+1)d}{2d(h+1)-4}, &\quad d\,{\rm even},\cr
     \frac{d+1}{2d-2}, &\quad d\,{\rm odd}, \cr}\right.   \lb{2.4}\ee
where $\dim M=hd$, $h=1$ when $M$ is a sphere $S^d$ of dimension $d$ and
\be \hat{\chi}(c) = \frac{1}{N(c)}
\sum_{0\le l_m\le \nu(c^m) \atop 1\le m\le N(c)}(-1)^{i(c^m)+l_m}k_{l_m}^{\ep(c^m)}(c^m)\;
                 \in \;\Q.        \lb{2.5}\ee}

\subsection{The structure of $H_*(\Lm M/S^1, \Lm^0 M/S^1;\Q)$}

Set $\ol{\Lm}^0=\ol{\Lambda}^0M =\{{\rm
constant\;point\;curves\;in\;}M\}\cong M$. Let $(X,Y)$ be a
space pair such that the Betti numbers $b_i=b_i(X,Y)=\dim
H_i(X,Y;\Q)$ are finite for all $i\in \Z$. As usual the {\it
Poincar\'e series} of $(X,Y)$ is defined by the formal power series
$P(X, Y)=\sum_{i=0}^{\infty}b_it^i$. We need the following well
known version of results on Betti numbers.

{\bf Lemma 2.5.} (cf. Theorem 2.4 and Remark 2.5 of \cite{Rad1}, cf.
also Proposition 2.4 of \cite{LoD1}) {\it Let $(S^d,F)$ be a
$d$-dimensional Finsler sphere.}

(i) {\it When $d$ is odd, the Betti numbers are given by
\bea b_j
&=& \rank H_j(\Lm S^d/S^1,\Lm^0 S^d/S^1;\Q)  \nn\\
&=& \left\{\matrix{
    2,&\quad {\it if}\quad j\in \K\equiv \{k(d-1)\,|\,2\le k\in\N\},  \cr
    1,&\quad {\it if}\quad j\in \{d-1+2k\,|\,k\in\N_0\}\bs\K,  \cr
    0 &\quad {\it otherwise}. \cr}\right. \lb{2.6}\eea
For any $k\in \N$ and $k\ge d-1$, there holds }
\bea \sum_{j=0}^k(-1)^jb_j
&=& \sum_{0\le 2j\le k}b_{2j}   \nn\\
&=& [\frac{k}{d-1}] + [\frac{k}{2}] - \frac{d-1}{2}  \nn\\
&=& \frac{k(d+1)}{2(d-1)} - \frac{d-1}{2} - \ep_{d,1}(k)  \nn\\
&\le& \frac{k(d+1)}{2(d-1)} - \frac{d-1}{2}.   \lb{2.7}\eea
where $\ep_{d,1}(k) = \{\frac{k}{d-1}\} + \{\frac{k}{2}\}\in [0,\frac{3}{2}-\frac{1}{2(d-1)})$.

(ii) {\it When $d$ is even, the Betti numbers are given by
\bea b_j
&=& \rank H_j(\Lm S^d/S^1,\Lm^0 S^d/S^1;\Q)  \nn\\
&=& \left\{\matrix{
    2,&\quad {\it if}\quad j\in \K\equiv \{k(d-1)\,|\,3\le k\in (2\N+1)\},  \cr
    1,&\quad {\it if}\quad j\in \{d-1+2k\,|\,k\in\N_0\}\bs\K,  \cr
    0 &\quad {\it otherwise}. \cr}\right.    \lb{2.8}\eea
For any $k\in \N$ and $k\ge d-1$, there holds
\bea -\sum_{j=0}^k(-1)^jb_j
&=& \sum_{0\le 2j-1\le k}b_{2j-1}   \nn\\
&=& \left[\frac{[\frac{k}{d-1}]+1}{2}\right] + [\frac{k+1}{2}] - \frac{d}{2}  \nn\\
&=& \frac{k d}{2(d-1)} - \frac{d-2}{2} - \{\frac{[\frac{k}{d-1}]+1}{2}\}
          - \{\frac{k+1}{2}\} - \frac{1}{2}\{\frac{k}{d-1}\}  \nn\\
&\le& \frac{k d}{2(d-1)} - \frac{d-2}{2}.   \lb{2.9}\eea}

\medskip

{\bf Proof.} It suffices to prove (\ref{2.7}) and (\ref{2.9}).

When $d$ is odd, for any $k\in\N$ and $m\in [0,d-1)$, we have
\bea \sum_{0\le j\le k(d-1)+m}b_{j}
&=& \sum_{0\le 2j\le k(d-1)+m}b_{2j}    \nn\\
&=& 2(k-1) +\frac{k(d-1)-(d-3)}{2} - (k-1) + [\frac{m}{2}] \nn\\
&=& k + \frac{(k-1)(d-1)}{2} + [\frac{m}{2}].  \nn\eea
Thus for any integer $k\ge d-1$, because $d$ is odd, we obtain
\bea \sum_{0\le 2j\le k}b_{2j}
&=& [\frac{k}{d-1}] + \frac{([\frac{k}{d-1}]-1)(d-1)}{2}
    + \left[\frac{k-[\frac{k}{d-1}](d-1)}{2}\right]  \nn\\
&=& [\frac{k}{d-1}] + [\frac{k}{2}] - \frac{d-1}{2}  \nn\\
&=& \frac{k(d+1)}{2(d-1)} - \frac{d-1}{2} - \{\frac{k}{d-1}\} - \{\frac{k}{2}\} \nn\\
&\le& \frac{k(d+1)}{2(d-1)} - \frac{d-1}{2}.  \nn\eea
This proves (\ref{2.7}).

When $d$ is even, for any odd $k\in\N$ and $m\in [0,2(d-1))$, we have
\bea \sum_{0\le j\le k(d-1)+m}b_j
&=& \sum_{0\le 2j-1\le k(d-1)+m}b_{2j-1}    \nn\\
&=& 2\,\frac{k-1}{2} + \frac{k(d-1)-(d-3)}{2} - \frac{k-1}{2} + [\frac{m}{2}] \nn\\
&=& \frac{k+1}{2} + \frac{(k-1)(d-1)}{2} + [\frac{m}{2}].  \nn\eea
Note that for an integer $l>0$ there holds
$$  2\left[\frac{l+1}{2}\right] - 1 = \left\{\matrix{
         l, &\quad {\rm for}\;\;l\in 2\N-1, \cr
         l-1, &\quad {\rm for}\;\;l\in 2\N. \cr}\right.  $$
Thus for any integer $k\ge d-1$, because $d$ is even, we obtain
\bea \sum_{0\le j\le k}b_j
&=& \sum_{0\le 2j-1\le k}b_{2j-1}    \nn\\
&=& \frac{1}{2}\left((2\left[\frac{[\frac{k}{d-1}]+1}{2}\right]-1)+1\right)
    + \frac{1}{2}\left((2\left[\frac{[\frac{k}{d-1}]+1}{2}\right]-1)-1\right)(d-1)   \nn\\
&&\qquad + \left[\frac{1}{2}\left(k-(2\left[\frac{[\frac{k}{d-1}]+1}{2}\right]-1)(d-1)\right)\right] \nn\\
&=& \left[\frac{[\frac{k}{d-1}]+1}{2}\right]
    + \left(\left[\frac{[\frac{k}{d-1}]+1}{2}\right]-1\right)(d-1)    \nn\\
&&\qquad + \left[\frac{k+(d-1)}{2}\right] - \left[\frac{[\frac{k}{d-1}]+1}{2}\right](d-1)  \nn\\
&=& \left[\frac{[\frac{k}{d-1}]+1}{2}\right] + \left[\frac{k+1}{2}\right] - \frac{d}{2}  \nn\\
&=& \frac{kd}{2(d-1)} - \frac{d-2}{2} - \{\frac{[\frac{k}{d-1}]+1}{2}\}
      - \{\frac{k+1}{2}\} - \frac{1}{2}\{\frac{k}{d-1}\}   \nn\\
&\le& \frac{kd}{2(d-1)} - \frac{d-2}{2}.  \nn\eea
This proves (\ref{2.9}). \hfill\hb

\medskip

For a compact and simply connected Finsler manifold $M$ with
$H^*(M;\Q)\cong T_{d,h+1}(x)$, when $d$ is odd, then $x^2=0$ and
$h=1$ in $T_{d,h+1}(x)$. Thus $M$ is rationally homotopy equivalent to
$S^d$ (cf. Remark 2.5 of \cite{Rad1} and \cite{Hin1}). Therefore, next we only
consider the case when $d$ is even.

Then we have the following result.

\medskip

{\bf Lemma 2.6.} (cf. Theorem 2.4 of \cite{Rad1}) {\it Let $M$ be a compact
simply connected manifold with $H^*(M;\Q)\cong T_{d,h+1}(x)$ for some integer
$h\ge 2$ and even integer $d\ge 2$. Let $D=d(h+1)-2$ and
\bea \Om(d,h) = \{k\in 2\N-1&\,|\,& iD\le k-(d-1)=iD+jd\le iD+(h-1)d\;  \nn\\
         && \mbox{for some}\;i\in\N\;\mbox{and}\;j\in [1,h-1]\}. \lb{2.10}\eea
Then the Betti numbers of the free loop space of $M$ defined by
$b_q = \rank H_q(\Lm M/S^1,\Lm^0 M/S^1;\Q)$ for $q\in\Z$ are given by
\be b_q = \left\{\matrix{
    0, & \quad \mbox{if}\ q\ \mbox{is even or}\ q\le d-2,  \cr
    [\frac{q-(d-1)}{d}]+1, & \quad \mbox{if}\ q\in 2\N-1\;\mbox{and}\;d-1\le q < d-1+(h-1)d, \cr
    h+1, & \quad \mbox{if}\ q\in \Om(d,h), \cr
    h, & \quad \mbox{otherwise}. \cr}\right.
\lb{2.11}\ee
For every integer $k\ge d-1+(h-1)d=hd-1$, we have
\bea \sum_{q=0}^kb_q
&=& \frac{h(h+1)d}{2D}(k-(d-1)) - \frac{h(h-1)d}{4} + 1 + \ep_{d,h}(k) \nn\\
&\le& h(\frac{D}{2}+1)\frac{k-(d-1)}{D} - \frac{h(h-1)d}{4} + 1 + \{\frac{D}{hd}\{\frac{k-(d-1)}{D}\}\}  \nn\\
&<&  h(\frac{D}{2}+1)\frac{k-(d-1)}{D} - \frac{h(h-1)d}{4} + 2,   \lb{2.12}\eea
where
\bea \ep_{d,h}(k)
&=& \{\frac{D}{hd}\{\frac{k-(d-1)}{D}\}\} - (\frac{2}{d}+\frac{d-2}{hd})\{\frac{k-(d-1)}{D}\}   \nn\\
&&\quad - h\{\frac{D}{2}\{\frac{k-(d-1)}{D}\}\} - \{\frac{D}{d}\{\frac{k-(d-1)}{D}\}\},  \lb{2.13}\eea
and there hold $\ep_{d,h}(k)\in (-(h+2),1)$ and $\ep_{d,1}(k)\in (-2,0]$ for all integer $k\ge d-1$. }

\medskip

{\bf Proof.} For a compact and simply connected Finsler manifold $M$ with
$H^*(M;\Q)\cong T_{d,h+1}(x)$ and some even integer $d$, the following Poincar\'{e} series was
computed out by Theorem 2.4 of \cite{Rad1}
\be \sum_{q=0}^{+\infty}b_q\,t^q \equiv P({\Lm}M/S^1,{\Lm}^0M/S^1)(t) =
t^{d-1}\left(\frac{1}{1-t^2}+\frac{t^{d(h+1)-2}}{1-t^{d(h+1)-2}}\right)\frac{1-t^{dh}}{1-t^d}.
    \lb{2.14}\ee
Thus we get
\bea \sum_{q=0}^{+\infty}b_q\,t^q
&=& t^{d-1}\left(\sum_{i=0}^{+\infty}t^{2i}+\sum_{i=1}^{+\infty}t^{iD}\right)
         \sum_{j=0}^{h-1}t^{jd}  \nn\\
&=& t^{d-1}\left(\sum_{j=0}^{h-1}\sum_{i=0}^{+\infty}t^{2i+jd}
          + \sum_{j=0}^{h-1}\sum_{i=1}^{+\infty}t^{iD+jd}\right).   \lb{2.15}\eea

For the first sum in (\ref{2.15}), we have
\bea \sum_{k\in\Z}u_k t^k
&\equiv& \sum_{j=0}^{h-1}\sum_{i=0}^{+\infty}t^{2i+jd}   \nn\\
&=& \sum_{j=0}^{h-2}(j+1)\sum_{2i=jd}^{(j+1)d-2}t^{2i} + h\sum_{i=0}^{+\infty}t^{(h-1)d+2i},
\lb{2.16}\eea
where (\ref{2.16}) is obtained by listing all items $t^{2i+jd}$ into a strip with $j$
running from $0$ to $h-1$ downwards and $i$ running from $0$ to $+\infty$ rightwards,
and then summing up all terms with the exponents $0$, $2$, $\ldots$, $d-2$, $d$,
$\ldots$, $(h-1)d-2$, $(h-1)d$, $\ldots$, respectively. Therefore we obtain
\be u_k = \left\{\matrix{
     0, &\quad \mbox{if}\;k\in 2\Z-1\;\;\mbox{or}\;\;k<0,  \cr
     [\frac{k}{d}]+1, &\quad \mbox{if}\;k\in 2\N_0\;\mbox{and}\;0\le k< (h-1)d, \cr
     h, &\quad \mbox{if}\;k\in 2\N_0\;\mbox{and}\;(h-1)d\le k.  }\right. \lb{2.17}\ee

For the second sum in (\ref{2.15}), because $d>1$, we have $D=d(h+1)-2>(h-1)d$. Thus
we have
\be  i D > (i-1)D+(h-1)d, \qquad \forall\; i\in \N.   \lb{2.18}\ee
Therefore every integer in $\Om(d,h)$ is covered precisely once by elements in
$\Om(d,h)$. Then let
\be \sum_{k\in\Z}v_k t^k \equiv \sum_{j=0}^{h-1}\sum_{i=1}^{+\infty}t^{iD+jd}
       = \sum_{i=1}^{+\infty}\sum_{j=0}^{h-1}t^{iD+jd}.  \lb{2.19}\ee
Here no any two terms in (\ref{2.19}) with different indices $(i,j)$ have the same
exponent. Thus we obtain
\be v_k = \left\{\matrix{
    1, &\quad \mbox{if}\;k\in iD+d\N\;\mbox{and}\; iD\le k \le iD + (h-1)d\;\;\mbox{for some}\;\;i\in\N, \cr
    0, &\quad \mbox{otherwise}. \cr}\right. \lb{2.20}\ee
Then from (\ref{2.15}), (\ref{2.16}) and (\ref{2.19}) we obtain
\be  b_q = u_{q-(d-1)} + v_{q-(d-1)}, \qquad \forall\; q\in\Z. \lb{2.21}\ee
together with (\ref{2.17}) and (\ref{2.20}), it yields (\ref{2.11}).

Because $D=d(h+1)-2>(h-1)d$, to get the sum (\ref{2.12}), by (\ref{2.21}) for any integers
$p\ge 1$ and $0\le m\le D-1$ we compute
\bea \sum_{j=0}^{pD+m}(u_j+v_j)
&=& \sum_{j=0}^{h-2}(j+1)\sum_{2i=jd}^{(j+1)d-2}1 + \left(h\sum_{(h-1)d\le 2i\le pD+m}1\right) \nn\\
&&\qquad + \sum_{i=1}^{p-1}\sum_{j=0}^{h-1}1 + 1 + [\frac{m}{d}] - [\frac{m}{hd}]     \nn\\
&=& \frac{h(h-1)}{2}\cdot\frac{d}{2} + h(\frac{pD-(h-1)d+2}{2}+[\frac{m}{2}])   \nn\\
&&\qquad + (p-1)h + 1 + [\frac{m}{d}] - [\frac{m}{hd}]  \nn\\
&=& h(\frac{D}{2}+1)p - \frac{h(h-1)d}{4} + 1 + h[\frac{m}{2}] + [\frac{m}{d}] - [\frac{m}{hd}], \lb{2.22}\eea
where on the right hand side of the first equality the first two sums come from $u_j$s, and
the third sum and the last three terms come from $v_j$s. The number $1$ there corresponds to
the term $t^{pD}$. Note that by the fact $0<m\le D-1$, we have $m<(h-1)d + 2d=(h+1)d$. But we
may have $m\ge hd$. If this happens, the term $[m/d]$ in the right hand side of the first
equality will be precisely one greater than it should be in (\ref{2.19}). Thus the term
$-[m/(hd)]$ is added to cancel this possible surplus $1$. Then $1+[\frac{m}{d}]-[\frac{m}{hd}]$
gives the total contribution of $v_j$s after the power $t^{(p-1)D+(h-1)d}$.

Therefore for every integer $k\ge hd-1$, letting
$$   p=[\frac{k-(d-1)}{D}] \quad \mbox{and} \quad m=k-(d-1)-pD,   $$
we obtain
\be  m = k-(d-1)-[\frac{k-(d-1)}{D}]D = \{\frac{k-(d-1)}{D}\}D < D. \lb{2.23}\ee
Then by (\ref{2.21}) we obtain
\be \sum_{q=0}^kb_q = \sum_{q=0}^k(u_{q-(d-1)} + v_{q-(d-1)}) = \sum_{j=-(d-1)}^{k-(d-1)}(u_j + v_j)
    = \sum_{j=0}^{k-(d-1)}(u_j + v_j),  \lb{2.24}\ee
where (\ref{2.17}), (\ref{2.20}) and the fact $d\ge 2$ are used.

Thus replacing $k-(d-1)=pD+m$ with the above $p$ and $m$ into (\ref{2.22}) and replacing
$[a]$ by $a-\{a\}$ for $a\in \R$ below, we obtain
\bea \sum_{q=0}^kb_q
&=& \sum_{j=0}^{pD+m}(u_j + v_j)   \nn\\
&=& h(\frac{D}{2}+1)[\frac{k-(d-1)}{D}] - \frac{h(h-1)d}{4}+ 1 + h[\frac{k-(d-1)-[\frac{k-(d-1)}{D}]D}{2}]    \nn\\
&&\qquad  + [\frac{k-(d-1)-[\frac{k-(d-1)}{D}]D}{d}] - [\frac{k-(d-1)-[\frac{k-(d-1)}{D}]D}{hd}]  \nn\\
&=& h(\frac{D}{2}+1)\frac{k-(d-1)}{D} - \frac{h(h-1)d}{4} + 1 + h\frac{k-(d-1)-[\frac{k-(d-1)}{D}]D}{2}    \nn\\
&&\qquad + \frac{k-(d-1)-[\frac{k-(d-1)}{D}]D}{d} - \frac{k-(d-1)-[\frac{k-(d-1)}{D}]D}{hd} \nn\\
&&\qquad - h(\frac{D}{2}+1)\{\frac{k-(d-1)}{D}\} - h\{\frac{k-(d-1)-[\frac{k-(d-1)}{D}]D}{2}\}  \nn\\
&&\qquad - \{\frac{k-(d-1)-[\frac{k-(d-1)}{D}]D}{d}\} + \{\frac{k-(d-1)-[\frac{k-(d-1)}{D}]D}{hd}\}  \nn\\
&=& h(\frac{D}{2}+1)\frac{k-(d-1)}{D} - \frac{h(h-1)d}{4} + 1 + \frac{hD}{2}\{\frac{k-(d-1)}{D}\}    \nn\\
&&\qquad + \frac{D}{d}\{\frac{k-(d-1)}{D}\} - \frac{D}{hd}\{\frac{k-(d-1)}{D}\}
         - h(\frac{D}{2}+1)\{\frac{k-(d-1)}{D}\}  \nn\\
&&\qquad - h\{\frac{D}{2}\{\frac{k-(d-1)}{D}\}\} - \{\frac{D}{d}\{\frac{k-(d-1)}{D}\}\}
         + \{\frac{D}{h d}\{\frac{k-(d-1)}{D}\}\}  \nn\\
&=& h(\frac{D}{2}+1)\frac{k-(d-1)}{D} - \frac{h(h-1)d}{4} + 1 + (\frac{D}{d}-h-\frac{D}{h d})\{\frac{k-(d-1)}{D}\}  \nn\\
&&\qquad + \{\frac{D}{hd}\{\frac{k-(d-1)}{D}\}\} - h\{\frac{D}{2}\{\frac{k-(d-1)}{D}\}\}
        - \{\frac{D}{d}\{\frac{k-(d-1)}{D}\}\}.     \lb{2.25}\eea
Then from
$$   \frac{D}{d}-h-\frac{D}{hd}=1-\frac{2}{d}-\frac{D}{hd} = -\frac{2}{d}-\frac{d-2}{hd},  $$
we obtain (\ref{2.12}). \hfill\hb

{\bf Remark 2.7.} When $d$ is even and $h=1$, the first equality of (\ref{2.12})
is exactly the third equality of (\ref{2.9}). In fact, in this case, there holds
$h=1$ and $D=2(d-1)$. So by (\ref{2.12})-(\ref{2.13}) we have
\bea
\sum_{q=0}^kb_q |_{h=1}
&=&\frac{h(h+1)d}{2D}(k-(d-1)) - \frac{h(h-1)d}{4} + 1+\ep_{d,1}(k) \nn\\
&=&\frac{k d}{2(d-1)} - \frac{d-2}{2} - \left(\{\frac{k-(d-1)}{2(d-1)}\}
     + \{\frac{k-(d-1)}{2}\}\right).  \lb{2.26}\eea
On the other hand, by (\ref{2.9}) and noting that $d$ is even, we have
\bea
&&\{\frac{[\frac{k}{d-1}]+1}{2}\}+\{\frac{k+1}{2}\}+\frac{1}{2}\{\frac{k}{d-1}\}\nn\\
&&\qquad =\{\{\frac{k+(d-1)}{2(d-1)}\}-\frac{1}{2}\{\frac{k}{d-1}\}\}
            +\frac{1}{2}\{\frac{k}{d-1}\}+\{\frac{k-(d-1)}{2}\}.    \lb{2.27}\eea
Note that no matter the integer $[\frac{k}{d-1}]$ is odd or even, we have always
$$ \{\frac{k+(d-1)}{2(d-1)}\} = \{\frac{1}{2}([\frac{k}{d-1}]+1)+\frac{1}{2}\{\frac{k}{d-1}\}\}
       \ge \frac{1}{2}\{\frac{k}{d-1}\}. $$
Thus (\ref{2.27}) yields
\be \{\frac{[\frac{k}{d-1}]+1}{2}\}+\{\frac{k+1}{2}\}+\frac{1}{2}\{\frac{k}{d-1}\}
    =\{\frac{k+(d-1)}{2(d-1)}\}+\{\frac{k-(d-1)}{2}\}.   \lb{2.28}\ee
Then (\ref{2.26}) and (\ref{2.28}) complete the proof of the above claim.

\setcounter{equation}{0}
\section{Morse indices of closed geodesics}

\subsection{Basic normal form decompositions of symplectic matrices
and precise index iteration formulae}

In \cite{Lon1} of 1999, Y. Long established the basic normal form
decomposition of symplectic matrices. Based on this result he
further established the precise iteration formulae of indices of
symplectic paths in \cite{Lon2} of 2000. These results form the
basis of our study on the Morse indices and homological properties
of closed geodesics. Here we briefly review these results:

As in \cite{Lon3}, denote by
\bea
N_1(\lm, a) &=& \left(\matrix{\lm & a\cr
                                0 & \lm\cr}\right), \qquad {\rm for\;}\lm=\pm 1, \; a\in\R, \lb{3.1}\\
H(b) &=& \left(\matrix{b & 0\cr
                      0 & b^{-1}\cr}\right), \qquad {\rm for\;}b\in\R\bs\{0, \pm 1\}, \lb{3.2}\\
R(\th) &=& \left(\matrix{\cos\th & -\sin\th \cr
                           \sin\th & \cos\th\cr}\right), \qquad {\rm for\;}\th\in (0,\pi)\cup (\pi,2\pi), \lb{3.3}\\
N_2(e^{\th\sqrt{-1}}, B) &=& \left(\matrix{ R(\th) & B \cr
                  0 & R(\th)\cr}\right), \qquad {\rm for\;}\th\in (0,\pi)\cup (\pi,2\pi)\;\; {\rm and}\; \nn\\
        && \qquad B=\left(\matrix{b_1 & b_2\cr
                                  b_3 & b_4\cr}\right)\; {\rm with}\; b_j\in\R, \;\;
                                         {\rm and}\;\; b_2\not= b_3. \lb{3.4}\eea
Here $N_2(e^{\th\sqrt{-1}}, B)$ is non-trivial if $(b_2-b_3)\sin\theta<0$, and trivial
if $(b_2-b_3)\sin\theta>0$ as defined in \cite{Lon2} and Definition 1.8.11 of \cite{Lon3}.
Note that symplectic paths with end matrices in these two cases have rather different index
iteration properties as proved in \cite{Lon2} (cf. Theorems 8.2.3 and 8.2.4 of \cite{Lon3}).
In \cite{Lon1}-\cite{Lon3}, all the matrices listed in (\ref{3.1})-(\ref{3.4}) are called
{\bf basic normal forms} of symplectic matrices.

As in \cite{Lon3}, given any two real matrices of the square block form
$$ M_1=\left(\matrix{A_1 & B_1\cr C_1 & D_1\cr}\right)_{2i\times 2i},\qquad
   M_2=\left(\matrix{A_2 & B_2\cr C_2 & D_2\cr}\right)_{2j\times 2j},$$
the $\diamond$-sum (direct sum) of $M_1$ and $M_2$ is defined by the
$2(i+j)\times2(i+j)$ matrix
$$ M_1\diamond M_2=\left(\matrix{A_1 & 0 & B_1 & 0 \cr
                                   0 & A_2 & 0& B_2\cr
                                   C_1 & 0 & D_1 & 0 \cr
                                   0 & C_2 & 0 & D_2}\right). $$

{\bf Definition 3.1.} (cf. \cite{Lon2} and \cite{Lon3}) {\it For every
$P\in\Sp(2d)$, the homotopy set $\Omega(P)$ of $P$ in $\Sp(2d)$ is defined by
$$ \Om(P)=\{N\in\Sp(2d)\,|\,\sg(N)\cap\U=\sg(P)\cap\U\equiv\Gamma\;\mbox{and}
                    \;\nu_{\om}(N)=\nu_{\om}(P)\, \forall\om\in\Gamma\}, $$
where $\sg(P)$ denotes the spectrum of $P$,
$\nu_{\om}(P)\equiv\dim_{\C}\ker_{\C}(P-\om I)$ for all $\om\in\U$.
The homotopy component $\Om^0(P)$ of $P$ in $\Sp(2d)$ is defined by
the path connected component of $\Om(P)$ containing $P$ (cf. p.38 of
\cite{Lon3}). }

Note that $\Om^0(P)$ defines an equivalent relation among symplectic
matrices. Specially two matrices $N$ and $P\in\Sp(2d)$ are homotopic
if $N\in\Om^0(P)$, and in this case we write $N\approx P$.

Then the following decomposition theorem is proved in \cite{Lon1}
and \cite{Lon2}

\medskip

{\bf Theorem 3.2.} (cf. Theorem 7.8 of \cite{Lon1}, Lemma 2.3.5 and
Theorem 1.8.10 of \cite{Lon3}) {\it For every $P\in\Sp(2d)$, there
exists a continuous path $f\in\Om^0(P)$ such that $f(0)=P$ and
\bea
f(1)
&=& N_1(1,1)^{\dm p_-}\,\dm\,I_{2p_0}\,\dm\,N_1(1,-1)^{\dm p_+}  \nn\\
 &&\dm\,N_1(-1,1)^{\dm q_-}\,\dm\,(-I_{2q_0})\,\dm\,N_1(-1,-1)^{\dm q_+} \nn\\
 &&\dm\,R(\th_1)\,\dm\,\cdots\,\dm\,R(\th_k)\,\dm\,R(\th_{k+1})\,\dm\,\cdots\,\dm\,R(\th_r)  \nn\\
 &&\dm\,N_2(e^{\aa_{1}\sqrt{-1}},A_{1})\,\dm\,\cdots\,\dm\,N_2(e^{\aa_{k_{\ast}}\sqrt{-1}},A_{k_{\ast}})  \nn\\
 &&\qquad\qquad \,\dm\,N_2(e^{\aa_{k_{\ast}+1}\sqrt{-1}},A_{k_{\ast}+1})\,\dm\,\cdots
  \,\dm\,N_2(e^{\aa_{r_{\ast}}\sqrt{-1}},A_{r_{\ast}}) \nn\\
 &&\dm\,N_2(e^{\bb_{1}\sqrt{-1}},B_{1})\,\dm\,\cdots\,\dm\,N_2(e^{\bb_{k_{0}}\sqrt{-1}},B_{k_{0}})  \nn\\
 &&\qquad\qquad \,\dm\,N_2(e^{\bb_{k_0+1}\sqrt{-1}},B_{k_0+1})\,\dm\,\cdots\,
              \dm\,N_2(e^{\bb_{r_{0}}\sqrt{-1}},B_{r_{0}}) \nn\\
 &&\dm\,H(2)^{\dm h_+}\,\dm\,H(-2)^{\dm h_-}, \lb{3.5}\eea
where $\frac{\th_{j}}{2\pi}\not\in\Q$ for $1\le j\le k$ and
$\frac{\th_{j}}{2\pi}\in\Q$ for $k+1\le j\le r$;
$N_2(e^{\aa_{j}\sqrt{-1}},A_{j})$'s are nontrivial basic normal
forms with $\frac{\aa_{j}}{2\pi}\not\in\Q$ for $1\le j\le k_{\ast}$
and $\frac{\aa_{j}}{2\pi}\in\Q$ for $k_{\ast}+1\le j\le r_{\ast}$;
$N_2(e^{\bb_{j}\sqrt{-1}},B_{j})$'s are trivial basic normal forms
with $\frac{\bb_{j}}{2\pi}\not\in\Q$ for $1\le j\le k_0$ and
$\frac{\bb_{j}}{2\pi}\in\Q$ for $k_0+1\le j\le r_0$; $p_-=p_-(P)$,
$p_0=p_0(P)$, $p_+=p_+(P)$, $q_-=q_-(P)$, $q_0=q_0(P)$,
$q_+=q_+(P)$, $r=r(P)$, $k=k(P)$, $r_{j}=r_{j}(P)$, $k_j=k_j(P)$ with $j=\ast, 0$ and
$h_+=h_+(P)$ are nonnegative integers, and $h_-=h_-(P)\in \{0,1\}$;
$\th_j$, $\aa_j$, $\bb_j \in (0,\pi)\cup (\pi,2\pi)$; these integers
and real numbers are uniquely determined by $P$ and satisfy}
\be p_- + p_0 + p_+ + q_- + q_0 + q_+ + r + 2r_{\ast} + 2r_0 + h_- + h_+ = d. \lb{3.6}\ee

\medskip

For $\tau>0$ and $d\in\N$ let
$$  \P_{\tau}(2d)=\{\ga\in C([0,\tau],\Sp(2d)\,|\,\ga(0)=I\}. $$
Based on Theorem 3.2, the homotopy invariance and symplectic
additivity of the indices, the following precise iteration formula
was proved in \cite{Lon2}:

\medskip

{\bf Theorem 3.3.} (cf. \cite{Lon2}, Theorem 8.3.1 and Corollary
8.3.2 of \cite{Lon3}) {\it Let $\ga\in\P_{\tau}(2d)$. Denote the
basic normal form decomposition of $P\equiv \ga(\tau)$ by
(\ref{3.5}). Then we have
\bea i(\ga^m)
&=& m(i(\ga)+p_-+p_0-r ) + 2\sum_{j=1}^rE\left(\frac{m\th_j}{2\pi}\right) - r   \nn\\
&&  - p_- - p_0 - {{1+(-1)^m}\over 2}(q_0+q_+) \nn\\
&& + 2\sum_{j=k_{\ast}+1}^{r_{\ast}}\vf\left(\frac{m\aa_j}{2\pi}\right) - 2(r_{\ast}-k_{\ast}), \lb{3.7}\\
\nu(\ga^m) &=& \nu(\ga) + {{1+(-1)^m}\over 2}(q_-+2q_0+q_+) + 2\vs(m,\ga(\tau)),    \lb{3.8}\\
\hat{i}(\ga) &=& i(\ga) + p_- + p_0 - r +
\sum_{j=1}^r\frac{\th_j}{\pi},   \lb{3.9}\eea where we denote by}
\bea
\vs(m,\ga(\tau)) &=& (r-k) - \sum_{j=k+1}^r\vf\left(\frac{m\th_j}{2\pi}\right)  \nn\\
&& + (r_{\ast}-k_{\ast}) -
\sum_{j=k_{\ast}+1}^{r_{\ast}}\vf\left(\frac{m\aa_j}{2\pi}\right)
             + (r_0-k_0) - \sum_{j=k_0+1}^{r_0}\vf\left(\frac{m\bb_j}{2\pi}\right).
             \lb{3.10}
\eea

By Theorems 8.1.4-8.1.7 and 8.2.1-8.2.4 on pp179-187 of \cite{Lon3}, we have specially

\medskip

{\bf Proposition 3.4.} {\it Every path $\ga\in\P_{\tau}(2)$ with end matrix being homotopic
to one of the following matrices must have odd index $i(\ga)$,
\be   N_1(1, b_1), \quad N_1(-1,b_2), \quad R(\th), \quad {\it or}\quad H(-2),   \lb{3.11}\ee
where $b_1=0$ or $1$, $b_2=0$ or $\pm1$, and $\th\in(0,\pi)\cup(\pi,2\pi)$. Paths
$\xi\in\P_{\tau}(2)$ with end matrix being homotopic to $N_1(1,-1)$ or $H(2)$, and $\eta\in\P_{\tau}(4)$
with end matrix being homotopic to $N_2(\omega,B)$ must have even indices $i(\xi)$ and $i(\eta)$.}

\medskip

{\bf Remark 3.5.} Note that all closed geodesics on a simply
connected manifold $M$ are orientable. Therefore, the Morse index of
a closed geodesic on $M$ equals the above Maslov-type index of a
symplectic path starting from identity $I$ and ending at
$P_c\in\Sp(2(\dim M-1))$ (cf. Theorem 1.1 of \cite{Liu1} and Theorem
3 of \cite{Wil1}). Next we will apply the precise iteration indices to
study properties of Morse indices of closed geodesics.

\subsection{The monotonicity of index growth}

In \cite{LoD1}, closed geodesics on Finsler manifold are classified
into two families, rational and irrational ones, as follows.

{\bf Definition 3.6.} (cf. Definitions 3.4 and 3.6 of \cite{LoD1}) {\it
A matrix $P\in\Sp(2d)$ is {\bf rational}, if no basic normal form in
(\ref{3.5}) of $P$ is of the form $R(\th)$ with $\th/\pi\in
\R\bs\Q$, and is {\bf irrational}, otherwise. Let
$\nu(P)=\dim_{\R}\ker_{\R}(P-I)$. $P$ is {\bf equally degenerate}, if
$\nu(P^m)=\nu(P)$ for all $m\in\N$. $P$ {\bf completely non-degenerate},
if $\nu(P^m)=0$ for all $m\ge 1$.

Let $(M,F)$ be a $d$-dimensional Finsler manifold. Let $c$ be an orientable
closed geodesic on $(M,F)$ whose linearized Poincar\'e map is
denoted by $P_c$ and then $P_c\in\Sp(2d-2)$. The closed geodesic $c$
is {\bf rational}, {\bf irrational}, {\bf equally degenerate}, or
{\bf completely non-degenerate}, if so is $P_c$. The {\bf analytical period}
$n(c)$ of $c$ is defined by
\be n(c) = \min\{k\in\N\,|\,\nu(c^k)=\max_{m\ge 1}\nu(c^m)\;\;{\it and}\;\;
                  i(c^{m+k})-i(c^{m})\in 2\Z, \;\;\forall\,m\in\N\}. \lb{3.12}\ee}
The following is also defined in \cite{LoD1} for any closed geodesic $c$ on $(M,F)$, let
\be  n_0(c) = \min\{k\in\N\,|\,\nu(c^k)=\max_{m\ge 1}\nu(c^m)\}.  \lb{3.13}\ee
We have the following result.

{\bf Lemma 3.7.} (cf. Lemma 3.5 of \cite{LoD1}) {\it Let $(M,F)$ be a $d$-dimensional
Finsler manifold. Let $c$ be an orientable closed geodesic on $M$ whose linearized Poincar\'e map is
denoted by $P_c$. There hold
\bea
&&  n(c) = n_0(c)\;\;{\it or}\;\;2n_0(c),  \lb{3.14}\\
&&  n(c) = 2n_0(c)\;\;{\it if\;and\;only\;if}\;\; q_-=0, \;\;\;h_-=1\;\;\;{\it and}\;\;\;
                 n_0(c)\;\;{\it is\;odd}, \lb{3.15}\eea
where $q_-=q_-(P_c)$ and $h_-=h_-(P_c)$ with $P_c$ defined in (\ref{3.5}). }

We need

{\bf Lemma 3.8.} {\it Let $(M,F)$ be a Finsler manifold and $c$ be
a prime orientable closed geodesic on $M$. Let $n=n(c)$ be the analytical period of $c$.
Suppose $m\in [1,n-1]$ satisfies

(i) $\nu(c) < \nu(c^m) < \nu(c^n)$, and

(ii) there exists no $k\in [1,m-1]$ satisfying $\;\;k|m$, $\;\;k|n\;$ and $\nu(c^k)=\nu(c^m)$.

\noindent Then $\;m|n\;$ must hold. }

{\bf Remark 3.9.} Lemma 3.8 is precisely Proposition 3.12 of \cite{LoD1} when
$c$ is rational and its orbit is isolated in closed geodesic orbits in $\Lm M$ in
addition. By carefully checking the proof of this Proposition 3.12, one can find
that it works also for irrational prime closed geodesics, and the condition
on isolatedness in closed geodesic orbits in $\Lm M$ is not necessary. Therefore
we omit the details of this proof here.

{\bf Lemma 3.10.} {\it Let $(M,F)$ be a compact Finsler manifold. Let $c$ be
an orientable closed geodesic on $M$ with analytical period $n=n(c)$. Then $n=n(c)$ is
precisely the integer $N$ in Lemma 2.3, i.e., there holds also
\be  K(c^{n+m})=K(c^m),\qquad \forall\,m\ge 1,   \lb{3.16}\ee}

{\bf Proof.} In fact, by the definition (\ref{2.3}) of $K(c^m)$, it suffices to prove
\be  k_j^{\ep(c^{nl+m})}(c^{nl+m})=k_j^{\ep(c^m)}(c^m),\qquad \forall\;j\in\Z,\;
        l\in \N_0,\;1\le m< n. \lb{3.17}\ee

Note firstly that by the definition of $n=n(c)$, for all $l\in\N_0$ and $1\le m<n$,
we have $i(c^{nl+m})-i(c^m)\in 2\Z$. It then yields
\be \ep(c^{nl+m})=(-1)^{i(c^{nl+m})-i(c)}=(-1)^{i(c^{m})-i(c)}=\ep(c^m). \lb{3.18}\ee
By the definition of $n=n(c)$, for these integers $l$ and $m$ we have also
\be \nu(c^{nl+m}) = \nu(c^m).  \lb{3.19}\ee
Therefore we need only to prove (\ref{3.17}) for $0\le j\le \nu(c^m)$.

By Lemma 3.8 we obtain some integer $p\in [1,m]$ such that both $p|m$, $p|n$, and
$\nu(c^p)=\nu(c^m)=\nu(c^{nl+m})$ hold. Then $p|(nl+m)$ holds and by (iii) of lemma
2.2, we obtain
\be k_j^{\ep(c^{nl+m})}(c^{nl+m})=k_j^{\ep(c^p)}(c^p)=k_j^{\ep(c^m)}(c^m),\qquad
         \forall\;j\in\Z.  \lb{3.20}\ee
The proof is complete. \hfill\hb

{\bf Definition 3.11.} {\it For every matrix $P\in\Sp(2d)$, using its
basic normal form decomposition (\ref{3.5}) we define
\be  \left\{\matrix{
\sg(P) = r+p_++p_0+q_-+q_0,   \cr
s(P) = r+p_-+p_0+q_++q_0+2(r_*-k_*). \cr}\right. \lb{3.21}\ee
Recall that we have defined in \cite{LoD1}:
\be p(P) = p_0(P)+p_-(P)+q_0(P)+q_+(P)+r(P)+2r_{\ast}(P).  \lb{3.22}\ee}

{\bf Lemma 3.12.} {\it Let $c$ be a closed geodesic with mean index $\hat{i}(c)>0$
on a compact simply connected Finsler manifold $(M,F)$ of dimension $d\ge 2$.
Denote the basic normal form
decomposition of the linearized Poincar\'e map $P_c$ of $c$ by (\ref{3.5}). Denote by
$n=n(c)$ the analytical period of $c$. Let $\sg(c)=\sg(P_c)$ given by Definition 3.11.
Then for any even integer multiple $T>0$ of $n$, we have}
\be  i(c^T) + \nu(c^n) = \sg(c) \qquad \mod\;\; 2.  \lb{3.23}\ee

{\bf Proof.} By Theorem 3.3, the definition of $n$, and the evenness of $T$,
we obtain
\bea
i(c^T) &=& T(i(c)+p_-+p_0-r) + 2\sum_{j=1}^rE(\frac{T\th_j}{2\pi})  \nn\\
       &&\qquad  -r -p_- -p_0 -q_0 - q_+ -2(r_{\ast}-k_{\ast}),   \lb{3.24}\\
\nu(c^n) &=& \nu(c) +q_- + 2q_0 +q_+ + 2\zeta(T,\ga(\tau)), \lb{3.25}\eea
where $\zeta(T,\ga(\tau))$ is given by (\ref{3.10}). Thus we obtain
\bea i(c^T) + \nu(c^n)
&=& \nu(c) -r -p_- -p_0 -q_- -q_0   \qquad (\mod \;2)  \nn\\
&=& p_- +2p_0 +p_+ -r -p_- -p_0 -q_- -q_0   \qquad (\mod \;2)  \nn\\
&=& \sg(c)  \qquad (\mod \;2).  \lb{3.26}\eea
This proves the lemma. \hfill\hb

\medskip

When the Morse index of a prime closed geodesic on a Finsler
manifold $M$ is not too small, Morse indices of iterations of this
closed geodesic satisfy the following monotonicity property.

\medskip

{\bf Theorem 3.13.} {\it Let $c$ be a closed geodesic on a compact
simply connected Finsler manifold $M$ of dimension $d\ge 2$ satisfying
\be i(c)+p_0+p_-\ge q_0+q_++r+2(r_*-k_*),  \lb{3.27}\ee
where we denote the basic normal form decomposition of the linearized Poincar\'e
map $P_c$ of $c$ by (\ref{3.5}). Then there holds $i(c^{m+1})\ge i(c^m)$ for all
$m\ge 1$. In particular, the condition (\ref{3.27}) holds if $i(c)\ge d-2$.}

\medskip

{\bf Proof.} By (\ref{3.7}) in Theorem 3.3, for any $m\ge1$, we have
\bea i(c^{m+1})-i(c^m)
&=& i(c)+p_-+p_0-r+2\sum_{j=1}^r\left[E\left(\frac{(m+1)\th_j}{2\pi}\right)
        - E\left(\frac{m\th_j}{2\pi}\right)\right]  \nn\\
&&\quad+\frac{(-1)^m-(-1)^{m+1}}{2}(q_0+q_+)
    + 2\sum_{j=k_{\ast}+1}^{r_{\ast}}\left[\vf\left(\frac{(m+1)\aa_j}{2\pi}\right)
      -\vf\left(\frac{m\aa_j}{2\pi}\right)\right]\nn\\
&\ge& i(c)+p_-+p_0-(q_0+q_++r+2(r_*-k_*)), \lb{3.28}\eea
which, together with the condition (\ref{3.27}), yields the desired claim.

On the other hand, by Proposition 3.4 and the homotopy invariance and
symplectic additivity of the index, we have
\be i(c)=p_-+p_0+q_-+q_0+q_++r+h_-\quad(\mod\,2). \lb{3.29}\ee
By (\ref{3.6}) with $d$ replaced by $d-1$, it yields
$q_0+q_++r+2(r_*-k_*)\le d-1$. If $q_0+q_++r+2(r_*-k_*)=d-1$,
there holds $p_-+p_0+p_++q_-+2k_*+2r_0+h_-+h_+=0$ by (\ref{3.6}), which implies that
$i(c)+d-1=0\ (\mod\,2)$ by (\ref{3.14}). Thus by $i(c)\ge d-2$, we obtain
\be  i(c)+p_-+p_0-(q_0+q_++r+2(r_*-k_*))\ge i(c)-(d-1)\ge 0.  \lb{3.30}\ee
This completes the proof of Theorem 3.13. \hfill\hb

\subsection{The quasi-monotonicity of index growth}

Note that the Morse indices of closed geodesics in general are not
monotone if the initial Morse index is small enough. For irrational closed
geodesics with enough irrational rotation terms, in this section we
establish a similar property, which we call {\it quasi-monotonicity}, to
replace the monotonicity of the indices.

For rational closed geodesics, the properties of Morse indices of
their iterations have been completely understood in \cite{LoD1}.
Here, we are interested in properties of Morse indices of iterations
of irrational closed geodesics. This needs properties of sequences of vectors
in $\R^n$ uniformly distributed mod one in number theory which can be found
in pages 5-6 of \cite{GrR1}

{\bf Definition 3.14.} (cf. pages 5-6 of \cite{GrR1}) {\it For given
$v=(v_1,\ldots, v_n)\in \R^n$, define $v$ mod $1$ to be the vector
$\{v\}=(\{v_1\}, \ldots, \{v_n\})$. The sequence of vectors $\{u_k\}_{k\in\N}$
with $u_k\in\R^n$ is {\bf uniformly distributed mod one} if for any
$0\le b_j<c_j<1$ for $j=1, 2, \ldots, n$, we have }
$$ \lim_{n\to \infty}\frac{1}{N}\,^{\#}\{k\le N\;|\;\{u_k\}\in \oplus [b_j,c_j)\}
     = \Pi_{j=1}^n(c_j-b_j). $$

{\bf Proposition 3.15} (Kronecker's result, cf. page 6 of \cite{GrR1})
{\it If $1, v_1, \ldots, v_n$ are linearly independent over $\Q$, then the vectors
$\{(kv_1, \ldots, kv_n)\}_{k\in\N}$ are uniformly distributed mod one on $[0,1]^n$.}

\medskip

For our purpose, we need the following definition.

\medskip

{\bf Definition 3.16.} {\it Let $v=(v_1,\ldots, v_n)\in (\R\bs\Q)^n$. For a vertex
$\chi\in \{0,1\}^n$ of $[0,1]^n$, we call $v$ uniformly distributed mod one near
$\chi$, if for any given $\ep\in (0,1/2)$, there exist infinitely many $m\in\N$
such that }
\be  |\{mv\}-\chi| < \ep.  \lb{3.31}\ee

Note that when some of $1, v_1, \ldots, v_n$ are linearly dependent
over $\Q$, the Proposition 3.15 does not hold in general. In this case,
the sequence $\{\{kv\}\,|\,k\ge 1\}$ in general is only uniformly distributed
on the intersections of some lower dimensional hyperplanes with $[0,1]^n$. For
example, let $v_1\in (0,1)\bs\Q$, and $v_2=1-v_1$. Then $1, v_1,
v_2$ are linearly dependent over $\Q$, and for $v=(v_1,v_2)$ the
sequence $\{\{kv\}\,|\,k\in\N\}$ is dense on the second diagonal
$\{(x,y)\in [0,1]^2\,|\,x+y=1\}$ of $[0,1]^2$. Specially $v$ is
uniformly distributed mod one near the vertexes $(0,1)$ and $(1,0)$,
but is not uniformly distributed mod one near the vertexes $(0,0)$
and $(1,1)$. For another extremal example: let $v_1=v_2\in (0,1)\bs\Q$.
Then $1, v_1, v_2$ are linearly dependent over $\Q$, and for $v=(v_1,v_2)$
the sequence $\{\{kv\}\,|\,k\in\N\}$ is dense on the diagonal
$\{(x,y)\in [0,1]^2\,|\,x=y\}$ of $[0,1]^2$. Specially $v$ is
uniformly distributed mod one near the vertexes $(0,0)$ and $(1,1)$,
but is not uniformly distributed mod one near the vertexes $(1,0)$
and $(0,1)$.

We need the following Theorem 11.1.2 of \cite{Lon3} (originally proved as
Theorem 4.2 of \cite{LoZ1}) to continue our study.

{\bf Proposition 3.17.} (Y. Long and C. Zhu \cite{LoZ1}) {\it Fix
$v=(v_1,\ldots,v_n)\in \R^n$. Let $H$ be the closure of the subset
$\{ \{mv\} \,|\, m\in\N\}$ in $\T^n$ and $V=T_0\pi^{-1}H$ be the
tangent space of $\pi^{-1}H$ at the origin in $\R^n$, where
$\pi:\R^n\to \T^n$ is the projection map. Define
$$  A(v) =V\bs\cup_{v_k\in\R\bs\Q}\{x=(x_1,\ldots,x_n)\in V\,|\, x_k=0\}. $$
Define $\psi(x)=0$ when $x\ge 0$ and $\psi(x)=1$ when $x<0$. Then for
any $a=(a_1,\ldots,a_n) \in A(v)$, the vector
$$  \chi = (\psi(a_1), \ldots, \psi(a_n))   $$
makes
$$  | \{Nv\}-\chi | < \ep. $$
holds for infinitely many $N\in\N$.

  Moreover, this set $A(v)$ possesses the following properties.

  (a) $A(v)\ne\emptyset$.

  (b) When $v\in \Q^n$, there holds $V=A(v)=\{0\}$.

  (c) When $v\in\R^n\bs\Q^n$, there hold $\dim V\geq 1$,
$0\not\in A(v)\subset V$, $A(v)=-A(v)$, and that $A(v)$ is open in $V$.

  (d) When $\dim V=1$, there holds $A(v)=V\bs\{0\}$.

  (e) When $\dim V\ge 2$, $A(v)$ is obtained from $V$ by deleting all the
coordinate hyperplanes with dimension strictly smaller than $\dim V$ from
$V$, and specially $\dim A(v)=\dim V$. }

Denote by $\hat{1}=(1, \ldots, 1)\in\R^n$. Define the opposite
vertex $\hat{\chi}$ of a vertex $\chi$ in $[0,1]^n$ by
\be \hat{\chi} = \hat{1} - \chi.  \lb{3.32}\ee
The following lemma is a generalization of the above example and will be useful
later.

{\bf Lemma 3.18.} {\it Let $v=(v_1,\ldots, v_n)\in (\R\bs\Q)^n$ be uniformly
distributed mod one near a vertex $\chi\in \{0,1\}^n$ of $[0,1]^n$. Then $v$ is
also uniformly distributed mod one near the opposite vertex $\hat{\chi}$ of
$\chi$. }

{\bf Proof.} For the given $v\in (\R\bs\Q)^n$, we apply Proposition 3.17 to
prove the lemma. Using notations in Proposition 3.17, we obtain
$A(v)\neq\emptyset$ by the conclusion (a) and the fact $v\in (\R\bs\Q)^n$.

Now using the function $\psi:\R\to \{0,1\}$ in Proposition 3.17, we
further define a map $\hat{\psi}$ from $A(v)$ to vertexes of $[0,1]^n$ by
$$  \hat{\psi}(a)=(\psi(a_1), \ldots, \psi(a_n)), \qquad \forall\;
             a=(a_1, \ldots, a_n)\in A(v). $$
Then we have the following two claims:

{\bf Claim (i)} {\it If $v$ is uniformly distributed mod one near a vertex $\chi$
of $[0,1]^n$, then there exists an $a\in A(v)$ such that $\chi=\hat{\psi}(a)$. }

In fact, let $H_0$ be the closure of the set $\{ \{mv\} \,|\, m\in\N\}$ in
$[0,1]^n$. It is well known that the closed set $H_0$ consists of only
finitely many connected components which are intersections of parallel
equal dimensional subspaces with $[0,1]^n$ (cf. descriptions in Sections
23.4 on page 508 and 23.10 on page 522 of \cite{HaW1}) and determined by the
integral linearly dependent relations satisfied by the irrational numbers
$\{v_1, \ldots, v_n\}$. Then we have $H=\pi(H_0)$ and $V$ in Proposition 3.17
can be identified as the linear subspace of $\R^n$ passing through $0$, parallel
to $H_0$, and satisfying $\dim V=\dim H_0$. Because the closed set $H_0$ consists
of only finitely many connected components, we can choose an $\ep\in (0,1/4)$
sufficiently small such that the ball $B_{\ep}(\chi)$ with radius $\ep$ centered
at the point $\chi$ in $\R^n$ has non-empty intersection with only one connected
component $H_1$ of $H_0$. Because $v$ is uniformly distributed mod one near the
vertex $\chi$, we can choose a sufficiently large $m\in \N$ such that
$\{mv\}\in B_{\ep}(\chi)\cap H_1$. Denote the point $\{mv\}$ by
$b=(b_1, \ldots, b_n)$. Then by the choice of $\ep$, we have either $b_i\in (0,1/4)$
or $b_i\in  (3/4,1)$ for each $i=1,\ldots,n$. Here the fact $v\in (\R\bs \Q)^n$ is
used. We define a new point $a=(a_1, \ldots, a_n)$ by
$$ a_i=\left\{\matrix{b_i, &\;\; {\rm if}\;\;b_i\in (0,1/4), \cr
                      b_i-1, &\;\; {\rm if}\;\;b_i\in (3/4,1), \cr}\right. $$
for all $i=1,\ldots, n$. Denote the segment connecting $\chi$ to $b$ by $l_1$, and
the straight line passing through $0$ and parallel to $l_1$ by $l_2$ (cf. The
definitions of $a$ and $b$ in Figure 3.1).

\begin{figure}[h]
\centering
  \setlength{\unitlength}{0.05 mm}%
  \begin{picture}(1113.7, 1449.5)(0,0)
  \put(0,0){\includegraphics{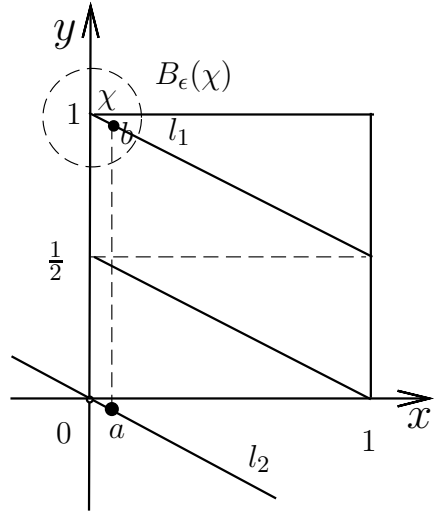}}
  \put(121.70,259.19){\fontsize{12.80}{15.36}\selectfont \makebox(45.0, 90.0)[l]{$0$\strut}}
  \put(230.63,1161.09){\fontsize{12.80}{15.36}\selectfont \makebox(45.0, 90.0)[l]{$\chi$\strut}}
  \put(1053.49,294.42){\fontsize{17.07}{20.48}\selectfont \makebox(60.0, 120.0)[l]{$x$\strut}}
  \put(112.99,1325.92){\fontsize{17.07}{20.48}\selectfont \makebox(60.0, 120.0)[l]{$y$\strut}}
  \put(418.73,1065.55){\fontsize{12.80}{15.36}\selectfont \makebox(90.0, 90.0)[l]{$l_1$\strut}}
  \put(631.48,193.84){\fontsize{12.80}{15.36}\selectfont \makebox(90.0, 90.0)[l]{$l_2$\strut}}
  \put(932.11,237.41){\fontsize{12.80}{15.36}\selectfont \makebox(45.0, 90.0)[l]{$1$\strut}}
  \put(147.85,1108.81){\fontsize{12.80}{15.36}\selectfont \makebox(45.0, 90.0)[l]{$1$\strut}}
  \put(291.63,1056.53){\fontsize{12.80}{15.36}\selectfont \makebox(45.0, 90.0)[l]{$b$\strut}}
  \put(261.13,272.26){\fontsize{12.80}{15.36}\selectfont \makebox(45.0, 90.0)[l]{$a$\strut}}
  \put(95.77,720.92){\fontsize{12.80}{15.36}\selectfont \makebox(135.0, 90.0)[l]{$\frac{1}{2}$\strut}}
  \put(383.12,1213.38){\fontsize{12.80}{15.36}\selectfont \makebox(45.0, 90.0)[l]{$B_\epsilon(\chi)$\strut}}
  \end{picture}%
\caption{\label{Fig3.1} The points $a$ and $b$ in the $2$ dimensional case.}
\end{figure}

Then by the definitions of $V$ and $A(v)$, we have
$$  a\;\in\; l_2\bs\{0\}\;\subset\; A(v)\;\subset\; V.  $$
Specially we have $\hat{\psi}(a)=\chi$ and proves the Claim (i).

{\bf Claim (ii)} {\it If a vertex $\chi$ of $[0,1]^n$ is in the image of $\hat{\psi}$,
so is its opposite vertex $\hat{\chi}$ in $[0,1]^n$. }

In fact, let $a=(a_1, \ldots, a_n)\in A(v)$ and $\chi=\hat{\psi}(a)$. Then $-a\in A(v)$ by
(c) of Proposition 3.17 and $a_i\not= 0$ for all
$1\le i\le n$ by (d) and (e) of Proposition 3.17. Therefore $\hat{\chi}=\hat{\psi}(-a)$ holds, which completes
the proof of Claim (ii).

Now Lemma 3.18 follows from these two claims. \hfill\hb

{\bf Corollary 3.19.} {\it Let $v=(v_1,\ldots, v_n)\in (\R\bs\Q)^n$. Then there exists
an integer $r$ satisfying $[(n+1)/2]\le r\le n$ and a subset $P$ of $\{1, \ldots, n\}$
containing $r$ integers, such that for any $\ep\in (0,1/4)$ there exist infinitely many
integers $T_1$ and $T_2\in\N$ satisfying respectively }
\bea
&& \left\{\matrix{ \{T_1v_i\} > 1-\ep, &\qquad& {\it for}\;\;i\in P, \cr
                     \{T_1v_j\} < \ep, &\qquad& {\it for}\;\; j\in \{1,\ldots,n\}\bs P, \cr}\right.  \lb{3.33}\\
&{\it and}\quad & \left\{\matrix{ \{T_2v_i\} < \ep, &\qquad& {\it for}\;\;i\in P, \cr
                     \{T_2v_j\} > 1- \ep, &\qquad& {\it for}\;\; j\in \{1,\ldots,n\}\bs P. \cr}\right.
                      \lb{3.34}\eea

{\bf Proof.} By the proof of Theorem 4.1 of \cite{LoZ1} (cf. pp.233-234 of \cite{Lon3}),
there exists a vertex $\chi$ of $[0,1]^n$ such that $v$ is uniformly distributed mod one near $\chi$.
By Lemma 3.18, $v$ is also uniformly distributed mod one near the opposite vertex $\hat{\chi}$ of
$\chi$ in $[0,1]^n$. Let $P(\xi)=\{j\in \{1, \ldots, n\}\,|\,\xi_j=1\}$ for any $\xi\in \{0,1\}^n$.
Let $P$ be the one of $P(\chi)$ and $P(\hat{\chi})$ which contains not fewer integers. Let
$r=\;^{\#}P$. Then the conclusion of Corollary 3.19 follows. \hfill\hb

\medskip

To estimate Morse indices of closed geodesics, we need first

{\bf Definition 3.20.} {\it For a prime orientable closed geodesic $c$ with mean index
$\hat{i}(c)>0$ on a Finsler manifold $(M,F)$ of dimension $d\ge 2$. Denote
the basic normal form decomposition of the linearized Poincar\'e map $P_c$
of $c$ by (\ref{3.5}). Using $\lm=i(c)+p_-+p_0-r$ and
$\rho(m)=\sum_{j=1}^r[\frac{m\th_j}{2\pi}]$ for any integer $A\in [0,k]$
with $k$ given in (\ref{3.5}), we define
\bea
\chi_c(m) &=& m\lm +2\rho(m), \qquad \forall\,m\in\N, \lb{3.35}\\
m_1(c) &=& \min\{\hat{m}\in\N\;|\;\chi_c(m)\ge i(c)+4\dim M + 2k,
                 \;\mbox{when}\ m\ge \hat{m}\}, \lb{3.36}\\
\aa_A(c) &=& \min\left\{\left.\left\{\frac{m\th_j}{2\pi}\right\}\,\right|\,
                1\le j\le A,\, 1\le m\le m_1(c)\right\},\lb{3.37}\\
\bb_A(c) &=& \min\left\{\left.\left\{\frac{m\th_j}{2\pi}\right\}\,\right|\,
                A+1\le j\le k,\, 1\le m\le m_1(c)\right\}.\lb{3.38}\eea}

Here we have $\aa_A(c)$ and $\bb_A(c)\in (0,1)$ whenever they are defined.
Note that from $\hat{i}(c)=\lm + \sum_{j=1}^r\frac{\th_j}{\pi}>0$, we obtain
\bea  \chi_c(m)
&=& m\lm + 2\sum_{j=1}^r[\frac{m\th_j}{2\pi}]   \nn\\
&=& m(\lm + \sum_{j=1}^r\frac{\th_j}{\pi}) - 2\sum_{j=1}^r\{\frac{m\th_j}{2\pi}\}   \nn\\
&=& m\hat{i}(c) - 2\sum_{j=1}^r\{\frac{m\th_j}{2\pi}\}   \nn\\
&\ge& m\hat{i}(c) - 2r.   \lb{3.39}\eea
Thus the positive integer $m_1(c)$ in (\ref{3.36}) and then $\aa_A(c)$ and $\bb_A(c)$
are well defined and depend only on $c$ , because $\hat{i}(c)>0$.

The following is our main estimate in this section.

\medskip

{\bf Theorem 3.21.} ({\bf Quasi-monotonicity of index growth for irrational closed geodesics})
{\it Let $c$ be a closed geodesic with mean index $\hat{i}(c)>0$ on a compact simply connected
Finsler manifold $(M,F)$
of dimension $d\ge 2$. Denote the basic normal form decomposition of the linearized Poincar\'e
map $P_c$ of $c$ by (\ref{3.5}). Then there exist an integer $A$ with $[(k+1)/2]\le A\le k$
and a subset $P$ of integers $\{1, \ldots, k\}$ with $A$ integers such that for any
$\ep\in (0,1/4)$ there exists an sufficiently large integer $T\in n\N$ satisfying
\bea
\left\{\frac{T\th_j}{2\pi}\right\} &>& 1-\ep, \qquad {\it for}\;\; j\in P,   \lb{3.40}\\
\left\{\frac{T\th_j}{2\pi}\right\} &<& \ep, \qquad {\it for}\;\; j\in \{1,\ldots,k\}\bs P.  \lb{3.41}\eea
Consequently we have
\bea
i(c^m)-i(c^T) &\ge& K_1 \equiv \lm + (q_0+q_+) + 2(r-k) + 2(r_{\ast}-k_{\ast}) + 2A,
          \quad \forall\, m\ge T+1, \lb{3.42}\\
i(c^T)-i(c^m) &\ge& K_2 \equiv \lm - (q_0+q_+) + 2k - 2(r_{\ast}-k_{\ast}) - 2A,
          \quad \forall\, 1\le m\le T-1, \lb{3.43}\eea
where $\lm = i(c)+p_-+p_0-r$, the integers $p_-$, $p_0$, $q_0$, $q_+$, $r$, $k$, $r_{\ast}$
and $k_{\ast}$ are defined in (\ref{3.5}). }

{\bf Remark 3.22.} Similar to our discussion after Definition 3.16, in Theorem 3.21 we
proved the existence of the integer $A$ located inside the interval $[[(k+1)/2],k]$. But
we do not know in general which precise value it may take without further knowledge on
the $\th_j/\pi$s. Specially if $k\ge 2$ and these irrational numbers are linearly dependent
over $\Q$, then $A$ can not take every integer value between $1$ and $k$.

\medskip

{\bf Proof of Theorem 3.21.} We carry out the proof in several steps.

{\bf Step 1.} Note first that iteration formulae of Morse indices of symplectic paths
ending at $N_1(-1,-1)$ or $N_1(-1,0)$ and those ending at $R(\pi)$ are precisely the same,
although their nullity may be different by $1$ (cf. Sections 8.1 and 8.2 of \cite{Lon3}).
Because our current theorem concerns only Morse indices of iterations of a closed geodesic,
so for simplicity of the description we shall replace all terms of $N_1(-1,-1)$ and
$N_1(-1,0)$ by $R(\pi)$s in (\ref{3.5}) and thus replace $r$ by the value of $r+q_0+q_+$
and set $\frac{1+(-1)^m}{2}(q_0+q_+)=0$ in (\ref{3.7}).

{\bf Step 2.} {\it Reduction to estimates on $\chi_c(m)$.}

Therefore, by (\ref{3.5}), Theorem 3.3 and Step 1, for any $m\ge 1$, the
iteration formulae of Morse indices of this closed geodesic $c$ is given by
\bea i(c^m)
&=& m(i(c)+p_-+p_0-r) + 2\sum_{j=1}^{r}E\left(\frac{m\th_j}{2\pi}\right)   \nn\\
&& - r - p_- - p_0 - 2(r_{\ast}-k_{\ast})
     + 2\sum_{j=k_{\ast}+1}^{r_{\ast}}\vf\left(\frac{m\aa_j}{2\pi}\right). \lb{3.44}\eea   
Therefore for any $T\in\N$, we obtain
\bea i(c^{m+T})
&-& i(c^T)\nn\\
&=& m\lm + 2\sum_{j=1}^{r}\left[E\left(\frac{(m+T)\th_j}{2\pi}\right)
               - E\left(\frac{T\th_j}{2\pi}\right)\right]
   + 2\sum_{j=k_{\ast}+1}^{r_{\ast}}\left[\vf\left(\frac{(m+T)\aa_j}{2\pi}\right)
                   - \vf\left(\frac{T\aa_j}{2\pi}\right)\right]  \nn\\
&=& m\lm + 2\left(\sum_{j=1}^{r}\left[\frac{m\th_j}{2\pi}\right]
     +\sum_{j=1}^{r}\left[E\left(\left\{\frac{T\th_j}{2\pi}\right\}+\left\{\frac{m\th_j}{2\pi}\right\}\right)
     -E\left(\left\{\frac{T\th_j}{2\pi}\right\}\right)\right]\right)   \nn\\
&{}&\qquad + 2\sum_{j=k_{\ast}+1}^{r_{\ast}}\left[\vf\left(\left\{\frac{T\aa_j}{2\pi}\right\}
                             +\left\{\frac{m\aa_j}{2\pi}\right\}\right)
                   - \vf\left(\left\{\frac{T\aa_j}{2\pi}\right\}\right)\right],  \quad
                   \forall\, m\ge 1, \lb{3.45}\eea
and
\bea i(c^{T})
&-& i(c^{T-m}) \nn\\
&=& m\lm +2\sum_{j=1}^{r}\left[E\left(\frac{T\th_j}{2\pi}\right)-E\left(\frac{(T-m)\th_j}{2\pi}\right)\right]
         +2\sum_{j=k_{\ast}+1}^{r_{\ast}}\left[\vf\left(\frac{T\aa_j}{2\pi}\right)
              -\vf\left(\frac{(T-m)\aa_j}{2\pi}\right)\right] \nn\\
&=& m\lm + 2\left(\sum_{j=1}^{r}\left[\frac{m\th_j}{2\pi}\right]
      + \sum_{j=1}^{r}\left[E\left(\left\{\frac{T\th_j}{2\pi}\right\}\right)
      - E\left(\left\{\frac{T\th_j}{2\pi}\right\}-\left\{\frac{m\th_j}{2\pi}\right\}\right)\right]\right) \nn\\
&{}&\qquad +
2\sum_{j=k_{\ast}+1}^{r_{\ast}}\left[\vf\left(\left\{\frac{T\aa_j}{2\pi}\right\}\right)
           - \vf\left(\left\{\frac{T\aa_j}{2\pi}\right\}-\left\{\frac{m\aa_j}{2\pi}\right\}\right)\right],
              \quad \forall\, 1\le m\le T-1.\lb{3.46}\eea

For $j=1,\ldots,r$, $i=1,\ldots,k_{\ast}$ and $m\ge 1$, let
\bea
\E_j^\pm(T,m) &=& E\left(\left\{\frac{T\th_j}{2\pi}\right\}
                  \pm\left\{\frac{m\th_j}{2\pi}\right\}\right), \lb{3.47}\\
\E_j(T) &=& E\left(\left\{\frac{T\th_j}{2\pi}\right\}\right), \lb{3.48}\\
\vf_i^\pm(T,m) &=& \vf\left(\left\{\frac{T\aa_i}{2\pi}\right\}
                   \pm\left\{\frac{m\aa_i}{2\pi}\right\}\right), \lb{3.49}\\
\vf_i(T) &=& \vf\left(\left\{\frac{T\aa_i}{2\pi}\right\}\right). \lb{3.50}\eea
Using these notations, from (\ref{3.45}) and (\ref{3.46}) we obtain
\bea i(c^{m+T}) - i(c^T)
&=& \chi_c(m) + 2\sum_{j=1}^{r}(\E_j^+(T,m) - \E_j(T))
    + 2\sum_{j=k_{\ast}+1}^{r_{\ast}}(\vf_j^+(T,m) - \vf_j(T)), \nn\\
&& \qquad\qquad \forall\, m\ge 1, \lb{3.51}\\
i(c^{T}) - i(c^{T-m})
&=& \chi_c(m) + 2\sum_{j=1}^{r}(\E_j(T) - \E_j^-(T,m))
    + 2\sum_{j=k_{\ast}+1}^{r_{\ast}}(\vf_j(T) - \vf_j^-(T,m)), \nn\\
&&\qquad\qquad \forall\, 1\le m\le T-1. \lb{3.52}\eea

Therefore (\ref{3.42}) and (\ref{3.43}) are equivalent to the following estimates:
\bea
\chi_c(m) &+& 2\sum_{j=1}^{r}(\E_j^+(T,m) - \E_j(T))
    + 2\sum_{j=k_{\ast}+1}^{r_{\ast}}(\vf_j^+(T,m) - \vf_j(T)) \ge K_1, \nn\\
&&\qquad\qquad \forall\, m\ge 1, \lb{3.53}\\
\chi_c(m) &+& 2\sum_{j=1}^{r}(\E_j(T) - \E_j^-(T,m))
    + 2\sum_{j=k_{\ast}+1}^{r_{\ast}}(\vf_j(T) - \vf_j^-(T,m)) \ge K_2, \nn\\
&&\qquad\qquad \forall\, 1\le m\le T-1. \lb{3.54}\eea

By the choice of $T\in n\N$, where $n=n(c)$ is the analytical period of $c$, we have
\bea
\E_j(T) &=& 1, \quad 1\le \E_j^+(T,m)\le 2, \quad 0\le \E_j^-(T,m)\le 1, {\hskip 1 cm} \forall\, 1\le j\le k,
                         \lb{3.55}\\
\E_j(T) &=& 0, \quad 0\le \E_j^+(T,m)\le 1, \quad \E_j^-(T,m)=0, {\hskip 1 cm} \forall\, k+1\le j\le r, \;\;m\ge 1,
                         \lb{3.56}\\
\vf_j(T) &=& 0, \qquad \vf_j^{\pm}(T,m) = \vf\left(\frac{m\aa_j}{2\pi}\right)\in \{0,1\},
   {\hskip 1 cm} \forall\, k_{\ast}+1\le j\le r_{\ast}, \;\;m\ge 1.   \lb{3.57}\eea
Therefore (\ref{3.53}) and (\ref{3.54}) are equivalent to the following estimates:
\bea
\chi_c(m) &+& 2\sum_{j=1}^{r}\E_j^+(T,m) - 2k
    + 2\sum_{j=k_{\ast}+1}^{r_{\ast}}\vf\left(\frac{m\aa_j}{2\pi}\right) \ge K_1, \nn\\
&&\qquad\qquad \forall\, m\ge 1, \lb{3.58}\\
\chi_c(m) &+& 2k - 2\sum_{j=1}^{k}\E_j^-(T,m)
    - 2\sum_{j=k_{\ast}+1}^{r_{\ast}}\vf\left(\frac{m\aa_j}{2\pi}\right) \ge K_2, \nn\\
&&\qquad\qquad \forall\, 1\le m\le T-1. \lb{3.59}\eea

{\bf Step 3.} {\it Definition of the set $S$.}

Note that by the Definition 3.20 of $m_1(c)$ and the definitions of $K_1$ and $K_2$ in
(\ref{3.42}) and (\ref{3.43}), we have
\be   \chi_c(m) \ge i(c) + 4\dim M + 2k \ge \max\{K_1, K_2+2(r_{\ast}-k_{\ast})\},
       \qquad\, \forall\,m\ge m_1(c). \lb{3.60}\ee
Thus together with (\ref{3.55})-(\ref{3.57}), the estimates (\ref{3.58}) and
(\ref{3.59}) hold for all $m\ge m_1(c)$. Therefore to continue the proof, it
suffices to find $T$ so that (\ref{3.58}) and (\ref{3.59}) hold for those $m\in \N$
satisfying
\be  1\le m\le m_1(c) \quad {\rm and}\quad
     \chi_c(m) < \max\{K_1, K_2+2(r_{\ast}-k_{\ast})\}.  \lb{3.61}\ee

According to (\ref{3.61}), let
\be  \hat{K} \equiv \max\{K_1, K_2+2(r_{\ast}-k_{\ast})\}-\lm.   \lb{3.62}\ee
Then
$$ \hat{K} = \max\{2(r-k+r_{\ast}-k_{\ast})+2A, \; 2(k-(r_{\ast}-k_{\ast}))-2A\} > 0, $$
where we have used the fact $q_0+q_+=0$ from Step 1 and the definitions of $K_1$ and $K_2$ in
(\ref{3.42}) and (\ref{3.43}).

For every integer $\mu \in [0,\hat{K}]$, let
\bea
\S_{\mu} &=& \left\{m\in\N\,\left|\ \chi_c(m) = \lm + \mu,\ 1\le m\le m_1(c)\right.\right\}, \lb{3.63}\\
\S &\equiv& \bigcup_{0\le \mu\le\hat{K}}\S_{\mu}.  \lb{3.64}\eea

For $m\in [1,m_1(c)]\bs S$, we have
\be  \chi_c(m) > \lm + \hat{K} \ge \max\{K_1, K_2+2(r_{\ast}-k_{\ast})\},
      \qquad \forall\;1\le m\le m_1(c). \lb{3.65}\ee
Therefore together with (\ref{3.55})-(\ref{3.57}), the estimates (\ref{3.58}) and (\ref{3.59}) hold
for all $m\in [1,m_1(c)]\bs S$.

Now it suffices to prove the estimates (\ref{3.58}) and (\ref{3.59}) for every $m\in S$.

{\bf Step 4.} {\it Determinations of $A$ and $T$. }

Now by Corollary 3.19, there exists an integer $A$ with $[(k+1)/2]\le A\le k$ and a subset $P$
of $\{1,\ldots, k\}$ with precisely $A$ integers satisfying the conditions (\ref{3.40}) and
(\ref{3.41}). Here specially the existence of the integer $T$ follows from the mod one uniformly
distribution property (cf. \cite{GrR1}) used in Corollary 3.19. For notational simplicity, by
reordering $\th_j$s, without loss of generality we assume $P=\{1, \ldots, A\}$ in the following.

Thus for $\aa_c(A)$ and $\bb_c(A)>0$ defined in (\ref{3.37}) and (\ref{3.38}), we can find
$T\in n\N$ such that
\bea
\left\{\frac{T\th_j}{2\pi}\right\} &>& 1-\aa_c(A), \qquad 1\le j\le A,  \lb{3.66}\\
\left\{\frac{T\th_j}{2\pi}\right\} &<& \bb_c(A), \qquad A+1\le j\le k.  \lb{3.67}\eea

Thus by (\ref{3.66}) and the definition (\ref{3.47}) of $\E_j^+(T,m)$ for
$1\le m\le m_1(c)$, we obtain
$$    1 < \{\frac{T\th_j}{2\pi}\}+\aa_c(A) < 2.  $$
Then for $1\le j\le A$ we have
\bea
\E_j^+(T,m)
&=& E\left(\left\{\frac{T\th_j}{2\pi}\right\}+\left\{\frac{m\th_j}{2\pi}\right\}\right)  \nn\\
&\ge& E\left(\left\{\frac{T\th_j}{2\pi}\right\}+\aa_c(A)\right)  \nn\\
&=& 2.     \lb{3.68}\eea
Therefore by (\ref{3.55}) for such a $j$ we obtain $\E_j^+(T,m)=2$, and then
\be \sum_{j=1}^r\E_j^+(T,m) \ge \sum_{j=1}^k\E_j^+(T,m)\ge 2A + (k-A) = k+A. \lb{3.69}\ee

Similarly by (\ref{3.67}) and the definition (\ref{3.47}) of $\E_j^-(T,m)$ for $1\le m\le m_1(c)$,
we obtain
$$  -1 < \left\{\frac{T\th_j}{2\pi}\right\}-\bb_c(A) < 0.  $$
Then for $A+1\le j\le k$ we have
\bea  \E_j^-(T,m)
&=& E\left(\left\{\frac{T\th_j}{2\pi}\right\}-\left\{\frac{m\th_j}{2\pi}\right\}\right)  \nn\\
&\le& E\left(\left\{\frac{T\th_j}{2\pi}\right\}-\bb_c(A)\right)  \nn\\
&=& 0. \lb{3.70}\eea
Therefore by (\ref{3.55}) for such a $j$ we obtain $\E_j^-(T,m)=0$, and then
\be \sum_{j=1}^k\E_j^-(T,m) \le \sum_{j=1}^A\E_j^-(T,m) \le A.  \lb{3.71}\ee

On the other hand, note that because $\frac{\th_j}{\pi}\not\in \Q$ holds for $j=1,\ldots,k$,
we obtain
\be \sum_{j=1}^kE\left(\left\{\frac{m\th_j}{2\pi}\right\}\right)=k,\qquad \forall\,m\ge 1. \lb{3.72}\ee
Note that
\bea i(c^m)
&=& m\lm + 2\sum_{j=1}^{r}E\left(\frac{m\th_j}{2\pi}\right) - r-p_--p_0 -2(r_{\ast}-k_{\ast})
   + 2\sum_{j=k_{\ast}+1}^{r_{\ast}}\vf\left(\frac{m\aa_j}{2\pi}\right)  \nn\\
&=& \chi_c(m) + 2\sum_{j=1}^{r}E\left(\{\frac{m\th_j}{2\pi}\}\right) - r-p_--p_0 -2(r_{\ast}-k_{\ast})
   + 2\sum_{j=k_{\ast}+1}^{r_{\ast}}\vf\left(\frac{m\aa_j}{2\pi}\right).   \lb{3.73}\eea
Thus by (\ref{3.72}) for every $m\in\N$ we obtain
\bea \chi_c(m)
&=& - 2\sum_{j=1}^{r}E\left(\left\{\frac{m\th_j}{2\pi}\right\}\right)
      + (r+p_-+p_0) + 2(r_{\ast}-k_{\ast})  \nn\\
&&\qquad  - 2\sum_{j=k_{\ast}+1}^{r_{\ast}}\vf\left(\frac{m\aa_j}{2\pi}\right) + i(c^m)  \nn\\
&=& \lm + (i(c^m)-i(c)) + 2(r-k-\sum_{j=k+1}^{r}E\left(\left\{\frac{m\th_j}{2\pi}\right\}\right)) \nn\\
&&\qquad  + 2(r_{\ast}-k_{\ast}-\sum_{j=k_{\ast}+1}^{r_{\ast}}\vf\left(\frac{m\aa_j}{2\pi}\right)).
         \lb{3.74}\eea

{\bf Step 5.} {\it Estimates (\ref{3.58}) and (\ref{3.59}) for $m\in S$. }

By the definition of $m\in S_{\mu}$ with $0 \le \mu \le \hat{K}$, we have $\chi_c(m) = \lm + \mu$.
Thus by (\ref{3.74}) for such an $m\in S_{\mu}$, (\ref{3.58}) and (\ref{3.59}) are equivalent to the following
estimates:
\bea
\sum_{j=1}^{r}\E_j^+(T,m)
&\ge& k - \frac{\lm}{2} - \frac{\mu}{2} - \sum_{j=k_{\ast}+1}^{r_{\ast}}\vf\left(\frac{m\aa_j}{2\pi}\right)
        + \frac{K_1}{2},     \forall\, 1\le m\le m_1(c), \lb{3.75}\\
\sum_{j=1}^{k}\E_j^-(T,m)
&\le& k + \frac{\lm}{2} + \frac{\mu}{2} - \sum_{j=k_{\ast}+1}^{r_{\ast}}\vf\left(\frac{m\aa_j}{2\pi}\right)
        - \frac{K_2}{2},     \forall\, 1\le m\le \min\{T-1,m_1(c)\}. \lb{3.76}\eea

We continue the study in three sub-steps according to the value of $\mu$ for (\ref{3.75}) and (\ref{3.76}).

{\bf Sub-step 1.} {\it Study on (\ref{3.75}) for $m\in S_{\mu}$ with $1\le \mu\le 2(r-k+r_{\ast}-k_{\ast})$. }

We start from the following

{\bf Claim 1:} {\it For any $m\in\S_{\mu}$ with $0\le \mu\le 2(r-k+r_{\ast}-k_{\ast})$, the set
\be \left\{\left\{\frac{m\th_j}{2\pi}\right\},\left.
           \left\{\frac{m\aa_l}{2\pi}\right\}\,\right|\,k+1\le j\le r,\,k_{\ast}+1\le l\le r_{\ast}\right\}
               \lb{3.77}\ee
contains at least $(r-k + r_{\ast}-k_{\ast})-[\mu/2]$ non-zero elements.}

In fact, if the claim does not hold, then the number of zero elements in $\S_{\mu}$ is at
least $[\mu/2]+1$.

By the Bott formula (cf. \cite{Bot1} and Section 12.1 of \cite{Lon3}), there always holds
$i(c^m)-i(c)\ge 0$ for all $m\in\N$. Thus by the definition of $\S_{\mu}$, (\ref{3.74}) and
the above assumption we obtain
\bea  \lm + \mu
&=& \chi_c(m) \nn\\
&=& \lm + 2\left(r - k - \sum_{j=k+1}^{r}E\left(\left\{\frac{m\th_j}{2\pi}\right\}\right)\right)
      + 2\left(r_{\ast}-k_{\ast} - \sum_{j=k_{\ast}+1}^{r_{\ast}}\vf\left(\frac{m\aa_j}{2\pi}\right)\right)
      + i(c^m)-i(c)    \nn\\
&\ge& \lm + 2([\mu/2]+1) + i(c^m)-i(c)  \nn\\
&\ge& \lm + \mu + 1 + i(c^m)-i(c).   \lb{3.78}\eea
This contradiction proves Claim 1.

Note that in this Sub-step 1, $\{\frac{T\th_j}{2\pi}\}=0$ for $k+1\le j\le r$ by the choice of $T\in n\N$.
Thus by Claim 1 we have
\bea  \sum_{j=k+1}^r\E_j^{+}(T,m)
&+& \sum_{j=k_{\ast}+1}^{r_{\ast}}\vf\left(\frac{m\aa_j}{2\pi}\right)   \nn\\
&=& \sum_{j=k+1}^{r}E(\{\frac{m\th_j}{2\pi}\})
       + \sum_{j=k_{\ast}+1}^{r_{\ast}}\vf\left(\frac{m\aa_j}{2\pi}\right)   \nn\\
&\ge& (r-k+r_{\ast}-k_{\ast}) - [\frac{\mu}{2}].   \lb{3.79}\eea
Therefore by (\ref{3.69}) and (\ref{3.79}), we obtain
\bea  \sum_{j=1}^{r}\E_j^{+}(T,m)
&=& \sum_{j=1}^{k}\E_j^{+}(T,m)
     + \left[\sum_{j=k+1}^{r}\E_j^{+}(T,m) + \sum_{j=k_{\ast}+1}^{r_{\ast}}\vf\left(\frac{m\aa_j}{2\pi}\right)\right]
     - \sum_{j=k_{\ast}+1}^{r_{\ast}}\vf\left(\frac{m\aa_j}{2\pi}\right)   \nn\\
&\ge& k + A + (r-k+r_{\ast}-k_{\ast}) - [\frac{\mu}{2}]
           - \sum_{j=k_{\ast}+1}^{r_{\ast}}\vf\left(\frac{m\aa_j}{2\pi}\right)\nn\\
&\ge& k - \frac{\lm}{2} - \frac{\mu}{2} - \sum_{j=k_{\ast}+1}^{r_{\ast}}\vf\left(\frac{m\aa_j}{2\pi}\right)
                + \frac{K_1}{2} \nn\\
&&\qquad\qquad  + \left(\frac{\lm}{2} + A + (r-k+r_{\ast}-k_{\ast}) - \frac{K_1}{2}\right).  \lb{3.80}\eea
Therefore to get (\ref{3.75}), we need to require the last line in the right hand side of
(\ref{3.80}) to be non-negative. Thus the largest value which $K_1$ can take to guarantee
(\ref{3.75}) is:
\be  K_1 = \lm + 2(r-k+r_{\ast}-k_{\ast}+A). \lb{3.81}\ee

{\bf Sub-step 2.} {\it Study on (\ref{3.75}) for $m\in S_{\mu}$ with $2(r-k+r_{\ast}-k_{\ast})< \mu\le \hat{K}$. }

In this case, by (\ref{3.69}) we have
\bea  \sum_{j=1}^{r}\E_j^{+}(T,m)
&=& \sum_{j=1}^{k}\E_j^{+}(T,m)
    + \left[\sum_{j=k+1}^{r}\E_j^{+}(T,m)+\sum_{j=k_{\ast}+1}^{r_{\ast}}\vf\left(\frac{m\aa_j}{2\pi}\right)\right]
    - \sum_{j=k_{\ast}+1}^{r_{\ast}}\vf\left(\frac{m\aa_j}{2\pi}\right)  \nn\\
&\ge& \sum_{j=1}^{k}\E_j^{+}(T,m) - \sum_{j=k_{\ast}+1}^{r_{\ast}}\vf\left(\frac{m\aa_j}{2\pi}\right)  \nn\\
&\ge& k + A - \sum_{j=k_{\ast}+1}^{r_{\ast}}\vf\left(\frac{m\aa_j}{2\pi}\right)  \nn\\
&=& k - \frac{\lm}{2} - \frac{\mu}{2} - \sum_{j=k_{\ast}+1}^{r_{\ast}}\vf\left(\frac{m\aa_j}{2\pi}\right)
       + \frac{K_1}{2}
    + \left(\frac{\lm}{2} + A + \frac{\mu}{2} - \frac{K_1}{2}\right).  \lb{3.82}\eea
Therefore the choice of $K_1$ by (\ref{3.81}) yields also the largest value which $K_1$
can take to guarantee $\frac{\lm}{2} + A + \frac{\mu}{2} - \frac{K_1}{2}\ge 0$ in the right hand
side of (\ref{3.82}), and then (\ref{3.75}).

{\bf Sub-step 3.} {\it (\ref{3.76}) for $m\in [1,\min\{T-1,m_1(c)\}]\cap\S_{\mu}$ with $0\le \mu\le \hat{K}$. }

In this case, for any $m\in [1,\min\{T-1,m_1(c)\}]\cap\S_{\mu}$ with $0\le \mu\le \hat{K}$, by
(\ref{3.55}) and (\ref{3.71}) we have
\bea \sum_{j=1}^k\E_j^{-}(T,m)
&\le& A  \nn\\
&=& k + \frac{\lm}{2} + \frac{\mu}{2}
       - \sum_{j=k_{\ast}+1}^{r_{\ast}}\vf\left(\frac{m\aa_j}{2\pi}\right)
                - \frac{K_2}{2} \nn\\
&& \qquad + A - k - \frac{\lm}{2} - \frac{\mu}{2}
       + \sum_{j=k_{\ast}+1}^{r_{\ast}}\vf\left(\frac{m\aa_j}{2\pi}\right)
             + \frac{K_2}{2} \nn\\
&\le& k + \frac{\lm}{2} + \frac{\mu}{2} -
\sum_{j=k_{\ast}+1}^{r_{\ast}}\vf\left(\frac{m\aa_j}{2\pi}\right)
            - \frac{K_2}{2} \nn\\
&& \qquad + \left(A - k - \frac{\lm}{2} - \frac{\mu}{2}
        + (r_{\ast}-k_{\ast}) + \frac{K_2}{2}\right). \lb{3.83}\eea
Therefore to make
$$  A - k - \frac{\lm}{2} - \frac{\mu}{2} + (r_{\ast}-k_{\ast}) + \frac{K_2}{2} \le 0, $$
the largest value which $K_2$ can take is
\be  K_2 = \lm + 2(k-A) - 2(r_{\ast}-k_{\ast}).  \lb{3.84}\ee
From (\ref{3.83}) we obtain that the choice (\ref{3.84}) of $K_2$ makes (\ref{3.76}) holds.

Now (\ref{3.81}) and (\ref{3.84}) make (\ref{3.75}) and (\ref{3.76}) hold, and complete
the proof.

{\bf Step 6.} As the final step, we come back to the discussion in the Step 1 of this proof,
i.e., we consider the quantity $q_0+q_+$ in (\ref{3.7}) and the constants $K_1$ and $K_2$.
Then replacing $r$ by $r+q_0+q_+$ in (\ref{3.81}) and (\ref{3.84}) we obtain
\bea
K_1 &=& (\lm -(q_0+q_+))+ 2(r+(q_0+q_+)-k+r_{\ast}-k_{\ast}+A)   \nn\\
    &=& \lm + (q_0+q_+) + 2(r-k+r_{\ast}-k_{\ast}+A),    \lb{3.85}\\
K_2 &=& \lm - (q_0+q_+) + 2(k-A) - 2(r_{\ast}-k_{\ast}). \lb{3.86}\eea
These two quantities yield (\ref{3.42}) and (\ref{3.43}) and complete the proof of
Theorem 3.21. \hfill\hb

\medskip

The following consequences of Theorem 3.21 will be used later in our proof.

\medskip

{\bf Corollary 3.23.} ({\bf Maximal index jump}) {\it Let $c$ be a closed
geodesic with mean index $\hat{i}(c)>0$ on a compact simply connected
Finsler manifold $(M,F)$ of dimension $d\ge 2$. Denote the basic normal form
decomposition of the linearized Poincar\'e map $P_c$ of $c$ by (\ref{3.5}).
Suppose the integer $A$ in Theorem 3.21 can be chosen to be equal to $k$
given by (\ref{3.5}). Then there exist infinitely many integers $T\in\N$ such that}
\bea i(c^{T+1})-i(c^T)
&=& \lm + (q_0+q_+) + 2r + 2(r_{\ast}-k_{\ast})  \nn\\
&=&  i(c) + p_- +p_0 + q_0 + q_+ + r + 2(r_{\ast}-k_{\ast}).  \lb{3.87}\eea

{\bf Proof.} If we change the equalities to inequalities so that the right
hand side of (\ref{3.87}) becomes a lower bound of $i(c^{T+1})-i(c^T)$, then
it follows from (\ref{3.42}) of Theorem 3.21 with $A=k$ immediately.

To get the equality, we choose $A=k$ in Definition 3.15 and (\ref{3.40}).
Together with (\ref{3.55})-(\ref{3.57}), with $T$ chosen by Theorem 3.21 we
obtain
\bea i(c^{T+1}) - i(c^T)
&=& \lm + 2\rho(1) + 2\sum_{j=1}^{r}(\E_j^+(T,1) - \E_j(T))  \nn\\
&&\qquad + 2\sum_{j=k_{\ast}+1}^{r_{\ast}}(\vf_j^+(T,1) - \vf_j(T)) + (q_0+q_+)   \nn\\
&=& \lm + 2k + 2(r-k) + 2(r_{\ast}-k_{\ast}) + (q_0+q_+), \nn\eea
which yields (\ref{3.87}) and completes the proof. \hfill\hb

\medskip

{\bf Corollary 3.24.} {\it Let $c$ be a completely non-degenerate closed geodesic
with mean index $\hat{i}(c)>0$ on a compact simply connected Finsler manifold
$(M,F)$ of dimension $d\ge 2$. Denote the basic normal form decomposition of the
linearized Poincar\'e map $P_c$ of $c$ by (\ref{3.5}). Let $r=k$ be the total number
of rotation matrices as in (\ref{3.5}). Then there exists an integer $A$ with
$[(r+1)/2]\le A\le r$ and infinitely many integers $T\in\N$ such that}
\bea
i(c^m)-i(c^T) &\ge& i(c) + (2A-r), \qquad\qquad\forall\,m\ge T+1,  \lb{3.88}\\
i(c^T)-i(c^m) &\ge& i(c) - (2A-r), \qquad\qquad\forall\,1\le m\le T-1.  \lb{3.89}\eea

\medskip

{\bf Proof.} Because $c$ is completely non-degenerate, we have
$p_-=p_0=p_+=q_-=q_0=q_+=0$, $\lm=i(c)-r$, $r=k$ and $r_{\ast}=k_{\ast}$. By
Theorem 3.21, we obtain (\ref{3.42}) and (\ref{3.43}) for some integer $A$ with
$[(r+1)/2]\le A\le r$ and some $T\in n\N$, where the constants $K_1$ and $K_2$ in
(\ref{3.42}) and (\ref{3.43}) are given by
\bea
K_1 &=& \lm + 2(r-k) + 2(r_{\ast}-k_{\ast}) + 2A = i(c) + (2A-r),  \lb{3.90}\\
K_2 &=& \lm + 2k - 2(r_{\ast}-k_{\ast}) - 2A = i(c) - (2A-r).  \lb{3.91}\eea
Therefore (\ref{3.42}) and (\ref{3.43}) yield (\ref{3.88}) and (\ref{3.89})
respectively. \hfill\hb

\medskip

{\bf Remark 3.25.} Note that Theorem 3.21 and all corollaries hold as well
for every symplectic path $\ga\in\P_{\tau}(2d)$ by our proofs above.
In addition, note that we can choose the $T$ to be some multiple of $n$ in
all above properties of Morse indices.

\setcounter{equation}{0}
\section{Rational and completely non-degenerate closed geodesics}

Let $(M,F)$ be a compact manifold with an
irreversible or reversible Finsler (including Riemannian) metric $F$. In this
section, we study closed geodesics on $M$. It is well known that if the total
number of prime closed geodesics on $M$ with a bumpy metric $F$ is finite, then every prime closed
geodesic $c$ must satisfy $\hat{i}(c)> 0$ by Theorem 2 of \cite{BaK1}. By
results of \cite{GrM1}, \cite{ViS1}, and Theorem 2.4 of \cite{Rad1}, we are
interested in compact simply connected manifolds. We start from some lemmas.

{\bf Lemma 4.1.} {\it Suppose that there exists only one prime closed geodesic
$c$ on a compact simply connected bumpy Finsler manifold $M$ with
$H^*(M;\Q)=T_{d,h+1}(x)$ for some integers $d\ge 2$ and $h\ge 1$. Then we have
\be  \hat{i}(c)>0, \quad i(c)=d-1 \quad \mbox{and}\quad M_q=b_q,
        \qquad \forall\; q\in\N_0.   \lb{4.1}\ee}

{\bf Proof.} It suffices to prove the last two claims in (\ref{4.1}).

If $i(c)+d$ is even, $(d+2j-1)-i(c)$ is odd. Then by Lemma 2.1 there holds
$\ol{C}_{d+2j-1}(E,c^m)=0$ for all $m\in\N$ and $j\in\Z$. And thus all Morse-type numbers
satisfy $M_{d+2j-1}=0$. But the Morse inequalities and Lemmas 2.5 and 2.6 then imply the
contradiction $0=M_{d-1}\ge b_{d-1}\ge 1$. So $i(c)+d$ must be odd.

Now by Lemma 2.1 again, we obtain $\ol{C}_{d+2j}(E,c^m)=0$ for all $m\in\N$ and
$j\in\Z$, because $(d+2j)-i(c)$ is odd. Thus $M_{d+2j}=0=b_{d+2j}$ by Lemmas 2.5 and 2.6
for all $j\in\Z$. Then $M_q=b_q$ for any $q\in\N_0$ follows from the Morse inequalities.

In addition, by Lemmas 2.5 and 2.6, it yields $b_{d-1}=1$ and $b_j=0$ for $j\le d-2$.
Thus by $M_{d-1}=b_{d-1}=1$ and Lemma 2.1, we get $i(c)=d-1$. \hfill\hb

\medskip

Note that in this section, we denote by $\Q^m$ the $m$ times of the module $\Q$
instead of using the notation $m\Q$ in order to make the text clearer.

\medskip

In this paper, when there is only one prime closed geodesic $c$ on a Finsler
manifold $(M,F)$, we denote the corresponding energy levels by $\ka_m=E(c^m)$
for $m\ge 1$.

{\bf Lemma 4.2.} (cf. Proposition 4.1 of \cite{LoD1}) {\it Let $(M,F)$ be a simply
connected compact Finsler manifold with $H^{\ast}(M,\Q)=T_{d,h+1}(x)$ and
possessing only one prime closed geodesic $c$ which is rational. Let $n=n(c)$ be
the analytical period of $c$. Denote by
$C_j=H_j(\ol{\Lm},\ol{\Lm}^{\ka_n})=\Q^{c_j}$ for all $j\in\Z$. Then there holds
\be  c_j = b_{j-i(c^n)-p(c)} \qquad \forall\;j\in\Z, \lb{4.2}\ee
where $b_j$ is the Betti numbers of the free loop space in Lemma 2.6, the constant $p(c)$
is defined by $p(c)=p(P_c)$ via the linearized Poncar\'e map $P_c$ of $c$ and definition
(\ref{3.22}). }

Next we give a slight modification of Theorem 5.2 of \cite{LoD1} to give a new result
which is designed for manifolds in above lemmas with some integer $h\ge 2$
and even integer $d\ge 2$. Here only the condition (\ref{4.7}) below is weakened
slightly than that in \cite{LoD1}.

{\bf Theorem 4.3.} {\it Let $(M,F)$ be a simply connected compact
Finsler manifold with $H^{\ast}(M,\Q)=T_{d,h+1}(x)$ and satisfying (OR) with the
only prime closed geodesic $c$. Let $n=n(c)$ be the analytical period of $c$.
Denote by
\be d_j = k_j^{\ep(c^n)}(c^n), \qquad \forall j\in\Z.   \lb{4.3}\ee
Suppose that there exist two integers $\mu\ge -1$ and $p(c)\ge 0$ such
that $c$ satisfies the following conditions:
\bea
&&  i(c^{m+n}) = i(c^n) + i(c^m) + p(c), \qquad \forall\; m \ge 1,  \lb{4.4}\\
&&  i(c^m)+\nu(c^m) \le i(c^n)+\mu,  \qquad \forall\; 1\le m<n, \lb{4.5}\\
&&  d_j = 0, \qquad \forall\;j\ge \mu+2,  \lb{4.6}\\
&&  H_{i(c^n)+\mu+1}(\ol{\Lm},\ol{\Lm}^{\ka_n})=0.     \lb{4.7}\eea
Then there exists an integer $\ka\ge 0$ such that
\be B(d,q)(i(c^n) + p(c)) + (-1)^{i(c^n)+\mu}\ka
     = \sum_{j=\mu-p(c)+1}^{i(c^n)+\mu}(-1)^j b_j. \lb{4.8}\ee}

{\bf Proof of Theorem 4.3.} We indicate necessary modifications of the proof of
Theorem 5.2 of \cite{LoD1} and are very sketchy here.

As in the Step 1 of the proof of Theorem 5.2 of \cite{LoD1}, for $j\in\Z$, we denote by
\be U_j=H_j(\ol{\Lm}^{\ka_n},\ol{\Lm}^0)=\Q^{u_j}, \quad B_j=H_j(\ol{\Lm},\ol{\Lm}^0)=\Q^{b_j},
          \quad C_j=H_j(\ol{\Lm},\ol{\Lm}^{\ka_n})=\Q^{c_j}.  \lb{4.9}\ee
Let $\bb=i(c^n)$. Then the long exact sequence of the triple
$(\ol{\Lm},\ol{\Lm}^{\ka_n},\ol{\Lm}^0)$ yields the following diagram:
\bea
\begin{tabular}{ccc ccc ccc ccc ccc}
$C_{\bb+\mu+1}$&$\to$&$U_{\bb+\mu}$&$\to$&$B_{\bb+\mu}$&$\to$&$C_{\bb+\mu}$
                     &$\to$&$\cdots$&$\to$&$U_0$&$\to$&$B_0$&$\to$&$C_0$ \\
$\parallel$&   &$\parallel$&   &$\parallel$&   &$\parallel$&  & & &$\parallel$& &$\parallel$& &$\parallel$ \\
$0$& &$\Q^{u_{\bb+\mu}}$&  &$\Q^{b_{\bb+\mu}}$& &$\Q^{c_{\bb+\mu}}$& &$\cdots$& &$\Q^{u_0}$& &$0$& &$\;0$, \\
\end{tabular}  \nn\eea
where $C_{\bb+\mu+1}=0=C_0$ follows from (\ref{4.7}) and Lemma 4.2, and $b_0=0$ follows from
Lemma 2.6. Then this long exact sequence yields
\be  0 = \sum_{j=0}^{\bb+\mu}(-1)^j(u_j - b_j + c_j).  \lb{4.10}\ee

Replacing (5.17) in \cite{LoD1} by the above (\ref{4.10}), repeating the proof
of Theorem 5.2 in \cite{LoD1} and using the above Lemma 4.2, we obtain
\bea  0
&=& B(d,q)(\bb + p(c)) - \sum_{j=0}^{\bb+\mu}(-1)^jb_j
      + \sum_{j=0}^{\bb+\mu}(-1)^jb_{j-\bb-p(c)} + (-1)^{\bb+\mu}u_{\bb+\mu+1}   \nn\\
&=& B(d,q)(\bb + p(c)) - \sum_{j=\mu-p(c)+1}^{\bb+\mu}(-1)^jb_j
        + (-1)^{\bb+\mu}u_{\bb+\mu+1}.         \lb{4.11}\eea
That is, (\ref{4.8}) holds with $\ka=u_{\bb+\mu+1}\ge 0$. \hfill\hb

Our main result in this section generalizes the multiplicity results in \cite{LoD1} on
rational closed geodesics on spheres, and in \cite{DuL1} and \cite{Rad4} on bumpy spheres,
and \cite{Rad5} on bumpy $\CP^2$ to all compact simply connected manifolds.

{\bf Theorem 4.4.} {\it Let $M$ be a compact simply connected manifold with
$H^*(M;\Q)\cong T_{d,h+1}(x)$ for some integers $h\ge 1$ and $d\ge 2$. Let
$F$ be an irreversible Finsler metric on $M$ and $c$ be the only prime closed
geodesic on $M$. Then $c$ can be neither rational nor completely non-degenerate. }

{\bf Proof.} Note that when $d$ is odd, then $h=1$ by Remark 2.5 of
\cite{Rad1}. Note that when $h=1$, $M$ is rationally homotopic to the sphere $S^d$.
In this case the conclusion that $c$ can not be rational follows from \cite{LoD1},
the conclusion that $c$ can not be completely non-degenerate follows from \cite{DuL1}.
Therefore it suffices to prove the theorem for the integer $h\ge 2$ and even integer
$d\ge 2$. We continue the proof in two claims.

\medskip

{\bf Claim 1:} {\it $c$ is not rational.}

In fact, assuming that $c$ is rational, we follow ideas in the proof of Theorem 6.1
of \cite{LoD1} and prove the theorem by contradiction.

To generate the non-trivial $H_{d-1}(\Lm M/S^1,\Lm M^0/S^1;\Q)$
(cf. Lemmas 2.5 and 2.6), the prime closed geodesic $c$ must satisfy
\be  \hat{i}(c)>0, \quad 0\le i(c)\le d-1.  \lb{4.12}\ee
Let $n=n(c)$ be the analytical period of $c$.
By the periodicity property (A) of Theorem 3.7 of \cite{LoD1}, we have
\be  i(c^{mn}) = m\,i(c^n) + (m-1)p(c), \qquad \forall \; m\in\N. \lb{4.13}\ee
Thus by (\ref{4.12}) and Corollary 9.2.7 of \cite{Lon4} we have
$i(c^n) + p(c) = \hat{i}(c^n) = n \hat{i}(c)>0$. Note that $i(c^n)=p(c)$ $\mod\;2$
by (D) of Theorem 3.7 of \cite{LoD1}, thus we have
\be  i(c^n) + p(c) \in 2\N.  \lb{4.14}\ee

Let $\mu=p(c)+(dh-3)$. Then by (\ref{4.14}) we have
\be  i(c^n) + \mu \ge dh-1 \ge 3, \qquad i(c^n) + \mu \in 2\N-1.  \lb{4.15}\ee

Now we can verify the conditions (\ref{4.4})-(\ref{4.7}) of Theorem 4.3
as in the proof of Theorem 6.1 of \cite{LoD1}. Note that (\ref{4.4})-(\ref{4.5})
follow from Theorem 3.7 and Proposition 3.11 of \cite{LoD1}, (\ref{4.6}) follows from (B-2) of
Theorem 4.1 of \cite{LoD1}, and (\ref{4.7}) follows from Lemmas 2.6 and 4.2 and
the evenness of $\mu+1-p(c)$. Then by Theorem 4.3, we obtain for some integer $\ka\ge 0$:
\be B(d,h)(i(c^n) + p(c))  + (-1)^{i(c^n)+\mu}\ka
          = \sum_{j=\mu-p(c)+1}^{i(c^n)+\mu}(-1)^jb_j.  \lb{4.16}\ee

Thus by (\ref{4.15}) we obtain
\be  B(d,h)(i(c^n) + p(c)) \ge -\sum_{\mu-p(c)+1\le 2j-1\le i(c^n)+\mu}b_{2j-1}.
     \lb{4.17}\ee
By Lemma 2.4 we have
$$  B(d,h) = -\frac{h(h+1)d}{2D}<0.  $$
Thus from Theorem 3.7 of \cite{LoD1}, we have
\be    i(c^n)+\mu-(d-1)=i(c^n)+p(c)+dh-d-2 \in 2\N.  \lb{4.18}\ee
By (\ref{4.17}), (\ref{4.18}) and (\ref{2.12}) we obtain
\bea  i(c^n) + p(c)
&\le& -\frac{1}{B(d,h)}\sum_{\mu-p(c)+1\le 2j-1\le i(c^n)+\mu}b_{2j-1}    \nn\\
&=& \frac{2D}{h(h+1)d}\left(\sum_{0\le 2j-1\le i(c^n)+\mu}b_{2j-1}-\sum_{0\le 2j-1\le dh-2}b_{2j-1}\right).
             \lb{4.19}\eea

Letting $D=d(h+1)-2$. Note that because $i(c)+p(c)\ge 2$ by (\ref{4.14}), we have
\be  i(c^n)+\mu = i(c^n)+p(c)+dh-3 \ge d-1+(h-1)d. \lb{4.20}\ee
Thus by Lemma 2.6 we have
\be \sum_{0\le 2j-1\le i(c^n)+\mu}b_{2j-1} =
 \frac{h(h+1)d}{2D}(i(c^n)+\mu-(d-1)) - \frac{h(h-1)d}{4} + 1 +\ep_{d,h}(i(c^n)+\mu). \lb{4.21}\ee
On the other hand, because $dh-3< dh-1= d-1+(h-1)d$, by Lemma 2.6 we have
\bea  \sum_{0\le 2j-1\le dh-3}b_{2j-1}
&=& \sum_{d-1\le 2j-1\le dh-3}\left([\frac{2j-1-(d-1)}{d}]+1\right)    \nn\\
&=& \sum_{d\le 2j\le dh-2}[\frac{2j}{d}]   \nn\\
&=& \sum_{\frac{d}{2}\le j\le \frac{dh}{2}-1}[\frac{j}{d/2}]   \nn\\
&=& \frac{dh/2-1-(d/2-1)}{d/2}(\frac{dh/2-1-(d/2-1)}{d/2}+1)\frac{1}{2}\cdot\frac{d}{2} \nn\\
&=& \frac{dh(h-1)}{4}.  \lb{4.22}\eea
Therefore we get
\bea
&& \sum_{0\le 2j-1\le i(c^n)+\mu}b_{2j-1} - \sum_{0\le 2j-1\le dh-3}b_{2j-1}   \nn\\
&&\quad = \frac{h(h+1)d}{2D}(i(c^n)+\mu-(d-1)) - \frac{h(h-1)d}{4} + 1 +\ep_{d,h}(i(c^n)+\mu)  \nn\\
&&\quad\qquad\qquad - \sum_{d-1\le 2j-1\le dh-3}\left([\frac{2j-1-(d-1)}{d}]+1\right)  \nn\\
&&\quad = \frac{h(h+1)d}{2D}(i(c^n)+p(c)+dh-d-2) - \frac{dh(h-1)}{2} + 1 + \ep_{d,h}(i(c^n)+\mu). \lb{4.23}\eea

Then (\ref{4.19}) becomes
$$  i(c^n) + p(c)
\le i(c^n)+p(c) + dh -d-2 + \frac{2D}{h(h+1)d}\left(1-\frac{dh(h-1)}{2}+\ep_{d,h}(i(c^n)+\mu)\right), $$
that is,
\bea  \ep_{d,h}(i(c^n)+\mu)
&\ge& \frac{h(h+1)d}{2D}\left(d+2 + \frac{(h-1)D}{h+1} - dh - \frac{2D}{h(h+1)d}\right)   \nn\\
&=& \frac{dh-(d-2)}{dh+(d-2)}.  \lb{4.24}\eea

Note that by (\ref{4.18}) we have
\be i(c^n)+\mu-(d-1) = i(c^n)+p(c)+dh-d-2 = i(c^n)+p(c)-2d + D. \lb{4.25}\ee
Let $\eta\in [0,D/2-1]$ be an integer such that
\be \frac{2\eta}{D} = \{\frac{i(c^n)+p(c)-2d}{D}\} = \{\frac{i(c^n)+\mu-(d-1)}{D}\}. \lb{4.26}\ee

By the definition (\ref{2.13}) of $\ep_{d,h}(i(c^n)+\mu)$ and (\ref{4.26}), we obtain
\bea  \ep_{d,h}(i(c^n)+\mu)
&=& \{\frac{D}{dh}\{\frac{i(c^n)+\mu-(d-1)}{D}\}\}
         - (\frac{2}{d}+\frac{d-2}{dh})\{\frac{i(c^n)+\mu-(d-1)}{D}\}   \nn\\
&&\qquad - h\{\frac{D}{2}\{\frac{i(c^n)+\mu-(d-1)}{D}\}\}
         - \{\frac{D}{d}\{\frac{i(c^n)+\mu-(d-1)}{D}\}\}    \nn\\
&=& \{\frac{2\eta}{dh}\} - (\frac{2}{d}+\frac{d-2}{dh})\frac{2\eta}{D}
      - h\{\frac{2\eta}{2}\} - \{\frac{2\eta}{d}\}   \nn\\
&=& \{\frac{2\eta}{dh}\} - (\frac{2}{d}+\frac{d-2}{dh})\frac{2\eta}{D} - \{\frac{2\eta}{d}\}  \nn\\
&\equiv& \ep(2\eta).   \lb{4.27}\eea

Now we claim
\be  \ep(2\eta) < \frac{dh-(d-2)}{dh+(d-2)}, \qquad \forall\; 2\eta\in [0,dh-2].  \lb{4.28}\ee

In fact, we write
\be  2\eta = pd + 2m \qquad \mbox{with some}\;\; p\in \N_0, \;\; 2m\in [0,d-2]. \lb{4.29}\ee
Then from $pd+2m=2\eta \le dh-2=(h-1)d+d-2$ we have
\be  p\in [0, h-1].  \lb{4.30}\ee
Therefore in this case we obtain
\bea \ep(2\eta)
&=& \frac{pd+2m}{dh} - (\frac{2}{d}+\frac{d-2}{dh})\frac{pd+2m}{D} - \frac{2m}{d}  \nn\\
&=& \frac{p}{h} - \frac{(2h+d-2)p}{hD} + \frac{2m}{dh} - \frac{(2h+d-2)2m}{dhD} - \frac{2m}{d}  \nn\\
&=& \frac{p}{h}(1-\frac{2h+d-2}{D}) + \frac{2m}{d}(\frac{1}{h} - \frac{2h+d-2}{hD} - 1)  \nn\\
&=& \frac{p(d-2) - 2mh}{D}  \nn\\
&\le& \frac{(h-1)(d-2)}{D}.  \lb{4.31}\eea
Now if (\ref{4.28}) does not hold, we then obtain
$$ \frac{dh-(d-2)}{D} \le \ep(2\eta) \le \frac{(h-1)(d-2)}{D}, $$
that is,
$$  dh - d +2 \le dh - d + 2 - 2h. $$
Because $h\ge 2$, this yields a contradiction and completes the proof of (\ref{4.28}).

If $d=2$, then $\{\frac{2\eta}{d}\}=0$ holds in (\ref{4.28}). Thus the definition
(\ref{4.27}) implies
\be \ep(2\eta) \le \ep(dh-2), \qquad \forall\; 2\eta\in [dh,D-2].  \lb{4.32}\ee

If $d\ge 4$, for any $2\eta\in [dh,D-2]$, write $2\eta = pdh + 2m$ for some $p\in\N_0$ and
$2m\in [0,dh-2]$. Then from $D-2=(h+1)d-4=hd+d-4$ we obtain $p\le 1$ and $2m\le d-4$.
Thus we have
\bea \ep(2\eta)
&=& \frac{2m}{dh} - (\frac{2}{d}+\frac{d-2}{dh})\frac{pdh+2m}{D} - \frac{2m}{d}  \nn\\
&=& \ep(2m) - (\frac{2}{d}+\frac{d-2}{dh})\frac{pdh}{D}  \nn\\
&\le& \ep(2m).  \lb{4.33}\eea

Therefore from (\ref{4.28}), (\ref{4.32}) and (\ref{4.33}), we obtain that (\ref{4.28})
holds in fact for all integer $\eta\in [0,D/2-1]$. This contradicts (\ref{4.24})
and completes the proof of Claim 1.

\medskip

{\bf Claim 2:} {\it $c$ is not completely non-degenerate.}

In fact, assuming that $c$ is completely non-degenerate, in which case $(M,F)$ becomes
bumpy, we prove the theorem by contradiction.

Then by Theorems 3.2 and 3.3, we have the precise index iteration formulae
\be i(c^m) = m(i(c)-r)+2\sum_{j=1}^r\left[\frac{m\th_j}{2\pi}\right]+r,
      \qquad \mbox{where}\;\frac{\th_j}{2\pi}\in (0,1)\setminus\Q,\,1\le j\le r.  \lb{4.34}\ee

By Claim 1 and the mean index identity, we have $r\ge 2$. Note that Claim 2 was proved
in \cite{DuL1} and \cite{Rad4} when $d$ is odd or $h=1$, and in \cite{Rad5} when $d=h=2$.
Next we give the proof of Claim 2 in two cases for all the values of $d\ge 2$ and $h\ge 1$,
which yields also a new proof for the results in \cite{DuL1}, \cite{Rad4}, and \cite{Rad5}.

\medskip

{\bf Case 1:} {\it $H^*(M;\Q)=T_{d,h+1}(x)$ with $d=2$ and $h\ge 1$.}

In this case, by the index iteration formulae (\ref{4.34}) and Lemma 4.1, it yields
\bea
&& i(c)=d-1=1, \lb{4.35}\\
&& i(c^{2j})= i(c^2)\quad (\mod\,2),\quad i(c^{2j-1})= i(c)\quad (\mod\,2),\qquad\forall\; j\ge 1.
    \lb{4.36}\eea
By Lemma 2.6, for any odd integer $k\ge 2h+1$ the Betti numbers $b_j$ in this case satisfy
\bea \sum_{j=0}^kb_{j}
&=& h(h+1)\frac{k-1}{2h} - \frac{h(h-1)}{2} + 1 - (h+1)\{h\{\frac{k-1}{2h}\}\}   \nn\\
&=& \frac{(h+1)(k-h+1)}{2},  \lb{4.37}\eea
where we have used the fact $\{h\{\frac{2m}{2h}\}\}=0$ for any $m\in\Z$.

Note that, by Theorem 3.21 there exists an integer subset $P$ of $\{1, \ldots, r\}$ containing
$r_1$ integers with $[(r+1)/2]\le r_1\le r-1$, without loss of generality we assume
$P=\{1,\ldots, r_1\}$, such that for any given $\ep\in (0,1/4)$ there exists a sufficiently
large $T\in 2\N$ satisfying
\bea
1-\{\frac{T\th_j}{2\pi}\} < \frac{\ep}{r},  &\qquad& \forall\; 1\le j\le r_1, \lb{4.38}\\
\{\frac{T\th_j}{2\pi}\} < \frac{\ep}{r}, &\qquad& \forall\; r_1+1\le j\le r. \lb{4.39}\eea

Thus, by Lemma 4.1 and Corollary 3.24 with $A=r_1$, we can choose $T\in 2\N$ sufficiently
large such that $R\equiv i(c^{T})\ge 2h+1$ and
\bea
i(c^m)-i(c^T) &\ge& i(c)-r+2A = 1+2r_1-r,\qquad\qquad\forall\,m\ge T+1,  \lb{4.40}\\
i(c^T)-i(c^m) &\ge& i(c)+r-2A = 1-(2r_1-r),\quad\qquad\forall\,1\le m\le T-1. \lb{4.41}\eea

\medskip

{\bf Case 1-1:} {\it $R\equiv i(c^T)\in 2\Z+1$.}

In this subcase, because $T$ is even, $r$ must be odd by (\ref{4.34}). By Claim 1 and Lemma
2.4 we must have $r\ge 2$. Therefore together with (\ref{4.36}) we must have
\be  i(c^2)\in 2\Z+1, \quad \mbox{and}\quad 3\le r\in 2\N-1.  \lb{4.42}\ee
Here we have $n=n(c)=1$ in Lemmas 2.4 and 3.10. Then by the facts $B(2,h)=-\frac{h+1}{2}$ and
$\hat{i}(c)=1-r+\sum_{j=1}^r\frac{\th_j}{\pi}$ which follows from (\ref{4.34}), by Lemma 2.4
we get
\be 1-r+\sum_{j=1}^r\frac{\th_j}{\pi}=\frac{2}{h+1}.  \lb{4.43}\ee
Thus by (\ref{4.34}) we obtain
\bea R
&\equiv& i(c^T)  \nn\\
&=& T(1-r)+2\sum_{j=1}^r\left[\frac{T\th_j}{2\pi}\right]+r  \nn\\
&=& T(1-r)+2\sum_{j=1}^r\frac{T\th_j}{2\pi}-2\sum_{j=1}^r\left\{\frac{T\th_j}{2\pi}\right\}+r \nn\\
&<& T(1-r)+2\sum_{j=1}^r\frac{T\th_j}{2\pi}-(2r_1-r)+2\ep    \nn\\
&=& \frac{2T}{h+1}-(2r_1-r)+2\ep,  \lb{4.44}\eea
where the first inequality follows from (\ref{4.38}).

On the other hand, by (\ref{4.35}), (\ref{4.36}), (\ref{4.42}) and Lemma 2.1, every iteration
$c^m$ with $m\ge 1$ contributes $1$ to the corresponding Morse-type number $M_{i(c^m)}$. Let
\be  \tdR=R+2r_1-r-1.  \lb{4.45}\ee
Note that $\tdR\ge R$ holds and it is odd. Therefore, by (\ref{4.40}), (\ref{4.41}) and
Lemma 4.1, we have
\be  \sum_{j=0}^{\tdR} b_j = \sum_{j=0}^{\tdR}M_j = T.  \lb{4.46}\ee
By (\ref{4.37}) and (\ref{4.46}), we obtain
\be \frac{(h+1)(\tdR-h+1)}{2} = T.   \lb{4.47}\ee
Now combining (\ref{4.44}) and (\ref{4.47}) together, we get
\bea  R+(2r_1-r)
&<& \frac{2T}{h+1} + 2\ep    \nn\\
&=& \frac{2}{h+1}\cdot\frac{(h+1)(\tdR-h+1)}{2} + 2\ep    \nn\\
&=& \tdR-h+1+2\ep,   \lb{4.48}\eea
which implies $h < 1$. Contradiction!

\medskip

{\bf Case 1-2:} {\it $R\equiv i(c^T)\in 2\Z$.}

In this subcase, by (\ref{4.34}), Lemma 2.4, and Claim 1, similarly to the Case 1-1 we obtain
\be  i(c^2)\in 2\Z \quad \mbox{and} \quad 2\le r\in 2\N.  \lb{4.49}\ee
Thus $n=n(c)=2$ in Lemmas 2.4 and 3.10. Thus from $B(2,h)=-\frac{h+1}{2}$ and
$\hat{i}(c)=1-r+\sum_{j=1}^r\frac{\th_j}{\pi}$, similarly to (\ref{4.43}), by Lemma 2.4
we obtain
\be  1-r+\sum_{j=1}^r\frac{\th_j}{\pi} = \frac{1}{h+1}.  \lb{4.50}\ee
Thus similarly to the proof of (\ref{4.44}), we obtain
\be  R+2r_1-r < \frac{T}{h+1} + 2\ep.   \lb{4.51}\ee

On the other hand, it follows from (\ref{4.34}) and Lemma 2.1 that
every $c^{2m-1}$ with $m\ge 1$ contributes $1$ to the corresponding Morse-type
number $M_{i(c^{2m-1})}$ and every $c^{2m}$ with $m\ge1$ has no contribution
to any Morse-type numbers. Let
\be \tdR=R+2r_1-r.  \lb{4.52}\ee
Note that $\tdR\ge R$ holds and it is even. Therefore by
(\ref{4.40}), (\ref{4.41}), Lemmas 2.1, 2.6 and 4.1, we obtain
\be   \sum_{j=0}^{\tdR-1}b_j = \sum_{j=0}^{\tdR}b_j = \sum_{j=0}^{\tdR} M_j = \frac{T}{2}.
       \lb{4.53}\ee
Then by (\ref{4.37}) and (\ref{4.53}), we obtain
\be  (h+1)(\tdR-h)=T.  \lb{4.54}\ee

Now from (\ref{4.51}) and (\ref{4.54}) we obtain
\be R+2r_1-r < \frac{T}{h+1}+ 2\ep = \frac{(h+1)(\tdR-h)}{h+1} + 2\ep = \tdR-h + 2\ep, \lb{4.55}\ee
which implies $h<1$. Contradiction!

\medskip

{\bf Case 2:} {\it $H^*(M;\Q)=T_{d,h+1}(x)$ with even $d\ge 3$ and $h\ge 1$.}

\medskip

In this case, we have $i(c)=d-1$ by Lemma 4.1. By Corollary 3.24, as in the proof of
Case 1, we can choose sufficiently large $T\in 2n\N$ with $n=n(c)$ being the analytical
period of $c$ such that
\bea
R\equiv i(c^T) &\ge& 2(d-1) + d(h-1) + 1,  \lb{4.56}\\
i(c^m)-i(c^T) &\ge& d-1+2r_1-r, \qquad\qquad\forall\,m\ge T+1, \lb{4.57}\\
i(c^T)-i(c^m) &\ge& d-1-(2r_1-r), \qquad\qquad\forall\,1\le m\le T-1. \lb{4.58}\eea

Let $\tdR=R+2r_1-r-(d-1)$. Then it follows from (\ref{4.57}) and (\ref{4.58}) that
\bea
i(c^m)&\ge& \tdR+2(d-1) \ge \tdR+4,  \qquad\qquad\forall\,m\ge T+1,  \lb{4.59}\\
i(c^m) &\le& \tdR,  \qquad\qquad\forall\,1\le m\le T-1. \lb{4.60}\eea

If $d\ge 4$ and $h\ge 1$, by (\ref{4.59}) and (\ref{4.60}) we obtain
$$  \{\tdR+1, \ldots, \tdR+5\}\cap\{i(c^m)\,|\,m\ge 1\}=\{\tdR+1, \ldots, \tdR+5\}\cap\{i(c^T)\}.  $$
Here note that $R=i(c^T)$ may also miss all of $\tdR+1, \ldots, \tdR+5$.
Therefore every $c^m$ with $m\in \N\bs\{T\}$ has no contribution to the Morse-type
numbers $M_{\tdR+1}, \;\ldots,\;M_{\tdR+5}$. Note that by (\ref{4.56}), we have
$\tdR > d-1+d(h-1)$. Thus by Lemmas 2.5, 2.6 and 4.1 there holds
\be  2 \le \sum_{j=1}^5b_{\tdR+j} = \sum_{j=1}^5M_{\tdR+j} \le 1. \lb{4.61}\ee
Contradiction!

If $d=3$, we then have $h=1$. By (\ref{4.59}) and (\ref{4.60}) we obtain
$$ \{\tdR+1, \tdR+2, \tdR+3\}\cap\{i(c^m)\,|\,m\ge 1\}=\{\tdR+1, \tdR+2, \tdR+3\}\cap\{i(c^T)\}.  $$
Note that in this case, $R$ and $r$ have the same parity by the choice of $T$, and thus $\tdR$
is even. Similarly to (\ref{4.61}) by Lemmas 2.5 and 4.1 we then obtain
\be  2 = b_{\tdR+2} = \sum_{j=1}^3b_{\tdR+j} = \sum_{j=1}^3M_{\tdR+j} \le 1. \lb{4.62}\ee
Contradiction!

This completes the proof of Claim 2 and Theorem 4.4. \hfill\hb

\medskip

\setcounter{equation}{0}
\section{On $4$-dimensional compact simply connected irreversible Finsler manifolds}

In this section, we give the proof of the main Theorem 1.2 about closed geodesics
on $4$-dimensional compact simply connected irreversible Finsler manifolds.

By our discussion in Section 1 and Theorems A and B, it suffices to consider the
case of the $4$-dimensional compact simply connected manifold $M$ satisfying
$H^*(M;\Q)\cong T_{d,h+1}(x)$ for some integers $d\ge 2$ and $h\ge 1$ with $hd=4$.
Thus we consider only the following two cases:
\be d=4\;\;\mbox{and}\;\;h=1, \qquad \mbox{or}\qquad d=h=2.  \lb{5.1}\ee
In these two cases, by Lemma 2.4 we have correspondingly
\be B(4,1)= -\frac{2}{3}, \qquad B(2,2) = -\frac{3}{2}.  \lb{5.2}\ee

Suppose $F$ is an irreversible Finsler metric on $M$. Assume that $c$ is the only
prime closed geodesic on $(M,F)$, and we prove Theorem 1.2 by contradiction.

Denote the basic normal form decomposition of the linearized Poincar\'e
map $P_c$ of $c$ by (\ref{3.5}). By Theorem 4.4, the closed geodesic $c$ can be neither
rational nor completely non-degenerate. Because $\dim M=4$, together with Rademacher's
identity (Lemma 2.4), the basic normal form decomposition of $P_c$ must contain precisely
two rotation matrices $R(\th_j)$ with $\th_j/(2\pi)\in(0,1)\bs\Q$ for $j=1$ and $2$, and have
the following form:
\be  P_c\approx R(\th_1)\dm R(\th_2)\dm G,   \lb{5.3}\ee
where $G$ is one of the $2\times 2$ matrices listed below:
\be \left\{\matrix{
N_1(1,a) & \;\;\mbox{with} \;\; a=-1,\;0\;\mbox{or}\;1,  \cr
N_1(-1,b) & \;\;\mbox{with} \;\; b=-1,\;0\;\mbox{or}\;1,  \cr
R(\th_3) & \;\;\mbox{with}\;\; \frac{\th_3}{2\pi}\in ((0,1)\cap \Q)\bs \{\frac{1}{2}\}. \cr}\right.
         \lb{5.4}\ee
Specially in (\ref{3.5}) of $P_c$ we have
\be k=2. \lb{5.5}\ee
Note that by Lemma 2.4, the irrational numbers
\be   \sg_j=\frac{\th_j}{2\pi}   \lb{5.6}\ee
for $j=1$ and $2$ are always linearly dependent on $\Q$ in the following. The following
lemma studies the situation in more details.

{\bf Lemma 5.1.} {\it Suppose $\sg_j\in (0,1)\bs \Q$ for $j=1$ and $2$ satisfy
\be  \sg_1 + \sg_2 = \frac{q}{p},  \lb{5.7}\ee
for some $p, q\in\N$ and $(p,q)=1$. Then for any $m\in\N$ there holds
\be [m\sg_1] + [m\sg_2] = \left\{\matrix{
          [\frac{mq}{p}], & {\it if}\;\; \{m\sg_1\}<\{\frac{mq}{p}\}, \cr
          [\frac{mq}{p}]-1, & {\it if}\;\; \{m\sg_1\}>\{\frac{mq}{p}\}. \cr}\right. \lb{5.8}\ee
Specially there holds}
\be  [m\sg_1] + [m\sg_2] = [\frac{mq}{p}]-1, \qquad {\it when}\quad m\in p\N. \lb{5.9}\ee

{\bf Proof.} Note first that for $m\in\N$ we have
\bea  \{m\sg_1\} + \{m\sg_2\}
&=& \{m\sg_1\} + \{\frac{mq}{p} - m\sg_1\}   \nn\\
&=& \{m\sg_1\} + \{\{\frac{mq}{p}\} - \{m\sg_1\}\}   \nn\\
&=& \left\{\matrix{
    \{\frac{mq}{p}\}, & {\it if}\;\; \{m\sg_1\}<\{\frac{mq}{p}\}, \cr
    \{\frac{mq}{p}\}+1, & {\it if}\;\; \{m\sg_1\}>\{\frac{mq}{p}\}. \cr}\right. \lb{5.10}\eea
Thus we get
\bea   [m\sg_1] + [m\sg_2]
&=& m(\sg_1+\sg_2) - (\{m\sg_1\}+\{m\sg_2\})   \nn\\
&=& \frac{mq}{p} - (\{m\sg_1\}+\{m\sg_2\})   \nn\\
&=& [\frac{mq}{p}] + \{\frac{mq}{p}\} - (\{m\sg_1\}+\{m\sg_2\}).   \nn\eea
Together with (\ref{5.10}), it yields (\ref{5.8}).

When $m\in p\N$, we have always $\{m\sg_1\}>0=\{mq/p\}$ by the irrationality of $\sg_1$. This
completes the proof. \hfill\hb

We continue the proof of Theorem 1.2 in several steps according to the value of $i(c)$ and
the form of $G$.

\medskip

{\bf Step 1.} {\it $i(c)=0$.}

\medskip

By Theorems 8.1.4-8.1.7 of \cite{Lon3} (cf. Theorem 3.3 above), in this case we
must have $G=N_1(1,-1)$, because all the other choices of $G$ in (\ref{5.4}) yield
an odd $i(c)$ by Proposition 3.4. Thus by Theorem 3.3, we have the precise index formulae
\be  i(c^m) = -2m + 2([m\sg_1]+[m\sg_2]) + 2, \quad \mbox{and}\quad \nu(c^m)=1,
              \quad \forall\; m\in\N.  \lb{5.11}\ee
Then we have $i(c^m)\in 2\Z$ for all $m\ge 1$ and $n=n(c)=1$. Thus (\ref{5.2}) and Lemma 2.4 yield
\be -\frac{1}{|B(d,h)|}(k_0(c)-k_1^+(c)) = \hat{i}(c) = -2+2(\sg_1+\sg_2) > 0,  \lb{5.12}\ee
Then Lemma 2.2, (\ref{5.2}) and (\ref{5.12}) imply
\be  k_1^+(c^m)=k_1^+(c)=1,\quad k_0(c^m)=0, \quad \forall\; m\ge 1.  \lb{5.13}\ee

Therefore by (\ref{5.11}), (\ref{5.13}) and the Morse inequality, we obtain
$$ M_{2j}=0,\qquad b_{2j+1}=M_{2j+1}=\;^\#\{m\in\N:\,i(c^m)=2j\}, \quad \forall\;j\in\N_0.  $$
Specially we have $b_1=M_1>0$. Then we must have
\be   d=h=2, \qquad {\rm and}\qquad B(d,h)=-\frac{3}{2}.  \lb{5.14}\ee
Thus by Lemma 2.6 with $d=h=2$, we obtain
\be  M_{2j}=b_{2j}=0, \; M_1=b_1=1, \; M_3=b_3=2, \; M_{2j+5}=b_{2j+5}=3,
      \qquad \forall\, j\in\N_0. \lb{5.15}\ee

Next we estimate $i(c^m)$ using Lemma 5.1. From (\ref{5.12})-(\ref{5.14}) we obtain
\be   \sg_1 + \sg_2 = \frac{4}{3}. \lb{5.16}\ee
Then by Lemma 5.1 we obtain
\be  [\frac{4m}{3}]-1 \le [m\sg_1] + [m\sg_2] \le [\frac{4m}{3}], \qquad \forall\; m\in\N. \lb{5.17}\ee
Thus by (\ref{5.9}) and (\ref{5.11}) for $m=3k\in\N$ we obtain
\be  i(c^{3k}) = -6k + 2([3k\cdot \frac{4}{3}]-1) + 2 = 2k,  \quad \forall\;k\in\N.  \lb{5.18}\ee
By (\ref{5.11}) and (\ref{5.17}) for $m=3k+1\in\N$, we obtain
$$ -2(3k+1)+2([(3k+1)\frac{4}{3}]-1)+2 \le i(c^{3k+1}) \le -2(3k+1)+2[(3k+1)\frac{4}{3}]+2. $$
That is,
\be  2k \le i(c^{3k+1}) \le 2k+2, \qquad \forall\;k\in\N_0.  \lb{5.19}\ee
Similarly for $m=3k+2\in\N$, we obtain
$$ -2(3k+2)+2([(3k+2)\frac{4}{3}]-1)+2 \le i(c^{3k+2}) \le -2(3k+2)+2[(3k+2)\frac{4}{3}]+2. $$
It yields also
\be  2k \le i(c^{3k+2}) \le 2k+2, \qquad \forall\;k\in\N_0.  \lb{5.20}\ee

Now we have the following

{\bf Claim 1:} {\it Besides $i(c)=0$, there hold }
\bea
i(c^2) &=& i(c^3) = 2,   \lb{5.21}\\
i(c^{3m+1}) &=& i(c^{3m+2}) = i(c^{3m+3}) = 2m+2, \qquad \forall\; m\in\N. \lb{5.22}\eea

In fact, by $i(c)=0$, (\ref{5.11}), (\ref{5.13}) and $b_1=1$, we obtain
\be  i(c^m)\in 2\N,  \qquad \forall\;m\ge 2.  \lb{5.23}\ee
Thus by (\ref{5.23}), (\ref{5.20}), (\ref{5.18}), (\ref{5.13}) and (\ref{5.15}), we obtain
(\ref{5.21}).

Now by (\ref{5.13}) and (\ref{5.15}), from (\ref{5.18})-(\ref{5.20}) we obtain (\ref{5.22})
for $m=1$. Then by an induction argument on $m$ we get (\ref{5.22}) for all $m\in\N$ and
complete the proof of Claim 1.

Now from (\ref{5.11}), (\ref{5.16}) and (\ref{5.22}), for any $m\in\N$ we obtain
\bea 2m+2
&=& i(c^{3m+1})   \nn\\
&=& -2(3m+1) + 2([(3m+1)\sg_1]+[(3m+1)\sg_2]) + 2   \nn\\
&=& -6m + 2((3m+1)(\sg_1+\sg_2)) - 2(\{(3m+1)\sg_1\} + \{(3m+1)\sg_2\})  \nn\\
&=& 2m + \frac{8}{3} - 2(\{(3m+1)\sg_1\} + \{(3m+1)\sg_2\}). \nn\eea
That is
\be  \{(3m+1)\sg_1\} + \{(3m+1)\sg_2\} = \frac{1}{3}, \qquad \forall\; m\in\N. \lb{5.24}\ee
Similarly to the proof of Lemma 5.1, by (\ref{5.16}) again for all $m\in\N$ we obtain
\bea  \{(3m+1)\sg_1\} + \{(3m+1)\sg_2\}
&=& \{(3m+1)\sg_1\} + \{(3m+1)(\frac{4}{3}-\sg_1)\}    \nn\\
&=& \{(3m+1)\sg_1\} + \{\frac{1}{3}-(3m+1)\sg_1\}  \nn\\
&=& \{(3m+1)\sg_1\} + \{\frac{1}{3}-\{(3m+1)\sg_1\}\}.  \lb{5.25}\eea
Because $\sg_1$ is irrational, by a result of A. Granville and Z. Rudnick (cf. the final
remark on page 6. of \cite{GrR1}), the sequence $\{(3m+1)\sg_1\}$ for $m\in\N$ is uniformly
distributed mod one on $[0,1]$. Thus we can find some sufficiently large $m\in\N$ such that
\be  \frac{1}{3} <  \{(3m+1)\sg_1\} < 1. \lb{5.26}\ee
Plugging it into (\ref{5.25}) yields the following identity for this $m$:
\be  \{(3m+1)\sg_1\} + \{(3m+1)\sg_2\} = \{(3m+1)\sg_1\} + 1 + \frac{1}{3}-\{(3m+1)\sg_1\}
    = \frac{4}{3},   \lb{5.27}\ee
which contradicts (\ref{5.24}). This proves that the Case of $i(c)=0$ can not happen.

\medskip

{\bf Step 2.} {\it $i(c)=1$.}

In this case, by $i(c)=1$ and Proposition 3.4, the matrix $G$ in (\ref{5.3}) must be one
of the the following matrices:
\be   N_1(1,a), \quad N_1(-1,b), \quad {\it or}\quad R(\th_3),   \lb{5.28}\ee
where $a=0$ or $1$, $b=\pm 1$, and $\frac{\th_3}{2\pi}\in(0,1)\cap\Q$.

Next we continue our proof in three subcases according to the particular form of the matrix $G$.

\medskip

{\bf Case 2-1:} {\it $G=N_1(-1,-1)$ or $R(\th_3)$ with $\frac{\th_3}{2\pi}\in(0,1)\cap\Q$.}

In this case, the index iteration formulae of $G=N_1(-1,-1)$ and $R(\th_3)$ are the same,
and only their nullities are different. Thus we can use the index formula for
$G=R(\th_3)$ to cover all these two subcases. As before, we write $\sg_j=\th_j/(2\pi)$ for
$j=1, 2, 3$. Then by Theorem 3.3, for $m\ge 1$ we have
\bea
i(c^m) &=& -2m+2\sum_{j=1}^3E(m\sg_j)-3,  \lb{5.29}\\
\nu(c^m) &=& \frac{1+(-1)^m}{2}(q_++2q_0)+2(r-2)-2\vf(m\sg_3). \lb{5.30}\eea
Specially in this case, we have $n=n(c)\ge 2$ and
\be   \left\{\matrix{
&1 \le i(c^m)\in 2\N-1, \qquad \forall\; m\in\N, \cr
&\nu(c)=0, \cr
&\nu(c^m)=0, \qquad \mbox{if}\;\; \frac{m}{n}\notin\Z,  \cr
&\nu(c^{mn})= 2, \qquad \mbox{for}\;\; m\in\N,  \qquad \mbox{if}\;\; G=R(\th_3), \cr
&\nu(c^{mn})= 1, \qquad \mbox{for}\;\; m\in\N,  \qquad \mbox{if}\;\; G=N_1(-1,-1). \cr}\right.
     \lb{5.31}\ee
Thus we have
\be  k_0(c^m) = k_0(c)=1, \qquad \mbox{for}\;\;1\le m\le n-1. \lb{5.32}\ee

Next we distinguish two subcases of $d=4$ with $h=1$ and $d=h=2$.

\smallskip

{\bf Subcase 2-1-1.} {\it $d=4$ and $h=1$. }

\smallskip

In this case, the manifold is rationally homotopic to $S^4$. We have

\smallskip

{\bf Claim 2:} {\it $i(c^n)=1$, $k_1^+(c^{nm})=k_1^{+}(c^n)\equiv k_1 \ge 1$ and
$k_0(c^{nm})=k_0(c^n)=k_2^+(c^{nm})=k_2^+(c^n)=0$ for all $m\in\N$.}

\smallskip

In fact, assume $i(c^n)\ge 3$. Then $i(c^{mn})\ge i(c^n)\ge 3$ for all $m\ge 1$.
Together with $i(c)=1$ and $\nu(c)=0$, it yields that the Morse type numbers satisfy
$M_2=M_0=0$ and $M_1\ge 1$. Then by Lemma 2.5 with $d=4$ the Morse inequality yields a
contradiction
$-1\ge M_2-M_1+M_0\ge b_2-b_1+b_0=0$. So $i(c^n)=1$ must hold.

Assume $k_1^+(c^n)=0$, by $i(c^m)\in 2\N-1$ and $\nu(c^m)\le 2$, we obtain $M_0=0$,
$M_1\ge 1$ and
\be M_{2j}=\,^{\#}\{m\in\N\,|\,i(c^{mn})=2j-1\}\,k_1^+(c^n)=0, \qquad \forall\;j\ge 1.
             \lb{5.33}\ee
Then the Morse inequality yields a contradiction $-1\ge M_2-M_1+M_0\ge b_2-b_1+b_0=0$. So
$k_1^+(c^n)\ge 1$ must hold, and then $k_0(c^n)=k_2^+(c^n)=0$ by Lemma 2.2. Then by
Lemma 2.2 again we get Claim 2 for all $m\in\N$.

\medskip

In this case, for numbers in the basic normal form decomposition (\ref{3.5})
of $P_c$ we have $r_*=p_-=p_0=0$, $k=2$, $r+q_0+q_+=3$. Note that Lemma 2.4 yields
the linear dependency of $1, \sg_1, \sg_2$ over $\Q$. Therefore in Theorem 3.21
we must have $A=1$, and we can find sufficiently large $T\in n\N$ such that
\bea
R &\equiv& i(c^T)\ge 3,  \lb{5.34}\\
i(c^m) &\ge& R+2, \qquad \forall\,m\ge T+1,  \lb{5.35}\\
i(c^m) &\le& R, \qquad \forall\,1\le m\le T. \lb{5.36}\eea

Because all $i(c^m)$ are odd, by (\ref{5.35}), (\ref{5.36}), Claim 2 and
Lemma 2.1 we obtain that $c^m$s with $m\ge T+1$ have no contributions to
$M_j$s with $0\le j\le R+1$, and $c^m$s with $1\le m\le T$ have only
contributions to $M_j$s with $0\le j\le R+1$. More precisely,
$\sum_{2j-1=1}^RM_{2j-1}$ is completely contributed by $c^m$s with
all integer $m\le T$ which are not in $n\N$, and each $c^m$ contributes a $1$.
And $\sum_{2j=0}^{R+1}M_{2j}$ is completely contributed by $c^{mn}$s with
$m\in\N$ satisfying $mn\le T$, and each $c^{mn}$ contributes a $k_1$. Thus we
have
\bea  \sum_{j=0}^{R+1}(-1)^j M_j
&=& \sum_{2j=0}^{R+1}M_{2j}-\sum_{2j-1=1}^RM_{2j-1}   \nn\\
&=& \frac{T}{n}k_1-\left(T-\frac{T}{n}\right)   \nn\\
&=& \frac{k_1-(n-1)}{n}T.   \lb{5.37}\eea

On the other hand, by (\ref{5.2}), Claim 2 and Lemma 2.4, we obtain
\be  \frac{k_1-(n-1)}{n\hat{i}(c)} = -\frac{2}{3}.  \lb{5.38}\ee

Thus by (\ref{5.37}), (\ref{5.38}), the Morse inequality and Lemma 2.5, we obtain
\be -\frac{2T}{3}\hat{i}(c) = \sum_{j=0}^{R+1}(-1)^jM_j
       \ge \sum_{j=0}^{R+1}(-1)^jb_j = -\sum_{2k-1=1}^R b_{2k-1} \ge 1-\frac{2R}{3}.
           \lb{5.39}\ee
It implies
\be 2R-2T\hat{i}(c) \ge 3.  \lb{5.40}\ee
However, by (\ref{5.29}) we have
\bea 2R-2T\hat{i}(c)
&=& 2\left(-2T+2\sum_{j=1}^3E(T\sg_j)-3\right)
       - 2T\left(2\sum_{j=1}^3\sg_j-2\right)   \nn\\
&=& 4\sum_{j=1}^2(E(T\sg_j)-T\sg_j) - 6  \nn\\
&\le& 2.   \lb{5.41}\eea
It contradicts to (\ref{5.40}) and then completes the proof in this subcase.

\medskip

{\bf Subcase 2-1-2.} {\it $d=h=2$. }

In this case, the manifold is rationally homotopic to $\CP^2$. Note that Lemma 2.4
yields the linear dependency of $1, \sg_1, \sg_2$ over $\Q$. Therefore in Theorem
3.21 we must have $A=1$, and there exists some $T\in 3n\N$ such that the odd integer
$R\equiv i(c^T)$ satisfies $R\ge 5$ and we have
\bea
i(c^m) &\ge& R+2, \qquad \forall\,m\ge T+1,  \lb{5.42}\\
i(c^m) &\le& R, \qquad \qquad\forall\,1\le m\le T. \lb{5.43}\eea

Let $k_j\equiv k_j^+(c^n)$ for $j=0$, $1$ or $2$. The following Claim 3 is crucial.

{\bf Claim 3:} {\it $k_2^+(c^{nm})=k_2=0$ for all $m\in\N$.}

If $\nu(c^n)=1$, then Claim 3 holds automatically by Lemma 2.2.
Next we consider the case of $\nu(c^n)=2$.

Otherwise, we assume $k_2=1$. Then by Lemma 2.2 we have
\be  k_2^+(c^{nm})=k_2=1,\quad k_0(c^{nm})=k_0=k_1^+(c^{nm})=k_1=0,\quad \forall\; m\in\N.
              \lb{5.44}\ee
Then by (\ref{5.2}) the identity in Lemma 2.4 becomes
$$   \frac{-(n-1)-k_2}{n\hat{i}(c)} = B(2,2) = -\frac{3}{2}. $$
Thus by (\ref{5.29}) and $k_2=1$ it yields
\be  \sg_1 + \sg_2 + \sg_3 = \frac{4}{3}.   \lb{5.45}\ee

Since $i(c^m)\in 2\N-1$ by (\ref{5.29}), there holds $M_{2j}=0$ for all $j\in\N_0$
by (\ref{5.44}) and Lemma 2.2. Thus together with the Morse inequality, it yields
\be M_{2j}=0,\quad M_{2j+1}=b_{2j+1},\quad\forall\,j\in\N_0.  \lb{5.46}\ee
By (\ref{5.44}) and Lemma 2.2, we have
\be M_{2j-1}=\,^{\#}\{m\in\N\,|\, i(c^m)=2j-1, \nu(c^m)=0\}
            + \,^{\#}\{m\in\N\,|\, i(c^m)=2j-3, \nu(c^m)=2\}.\lb{5.47}\ee
Let $N_{R+2}=\,^{\#}\{m\in\N\,|\, i(c^m)=R, \nu(c^m)=2\}$. Then it follows from (\ref{5.46}),
(\ref{5.47}) and $b_{R+2}=3$ by Lemma 2.6 that
\be   N_{R+2}\le M_{R+2}=b_{R+2}=3.  \lb{5.48}\ee

It follows from (\ref{5.42})-(\ref{5.44}) and (\ref{5.47})-(\ref{5.48}) that
\be \sum_{j=0}^R M_j=\sum_{2j-1=1}^R M_{2j-1}=T-N_{R+2}\ge T-3. \lb{5.49}\ee

On the other hand, by Lemma 2.6 with $d=h=2$, specially (\ref{4.37}), we obtain
\be  \sum_{j=0}^R b_j = \sum_{2j-1=1}^R b_{2j-1} = \frac{3(R-1)}{2}.   \lb{5.50}\ee

So (\ref{5.44}) and (\ref{5.49})-(\ref{5.50}) yield
\be  \frac{3(R-1)}{2}\ge T-3.  \lb{5.51}\ee

By (\ref{5.29}) and the definition of $T\in n\N$ we obtain
$$  R=i(c^T) = -2T+2\sum_{j=1}^3E(T\sg_j)-3 = -2T+2\sum_{j=1}^3[T\sg_j]+1.  $$
Therefore by (\ref{5.45}) we get
\bea \frac{3(R-1)}{2}
&=& - 3T + 3\sum_{j=1}^3[T\sg_j]  \nn\\
&=& - 3T + 3T(\sg_1+\sg_2+\sg_3) - 3(\{T\sg_1\}+\{T\sg_2\})    \nn\\
&=& T - 3(\{T\sg_1\}+\{T\sg_2\})     \nn\\
&=& T-3.   \lb{5.52}\eea
Here the last equality follows from that $\{T\sg_1\}+\{T\sg_2\}\in (0,2)$, $R$ is
odd and $T$ is an integer multiple of $3$, and then $\{T\sg_1\}+\{T\sg_2\}$ must be
an integer and then is equal to $1$. Then (\ref{5.49}), (\ref{5.50}) and (\ref{5.52})
yield $N_{R+2}=3$. In other words, by (\ref{5.43}) and the definition of $N_{R+2}$
there exist two distinct integers $T_1$ and $T_2$
with $T_1<T_2<T$ such that $i(c^{T_1})=i(c^{T_2})=R$ and $\nu(c^{T_1})=\nu(c^{T_2})=2$.
Because $\sg_3\in (0,1)\cap\Q$, there holds $\sg_3=q/p$ with some $q<p\in \N$ with
$(p,q)=1$. Therefore $p\ge 2$ holds. Then we have $T-T_2\ge 2$ and $T_2-T_1\ge 2$, thus
\be  T-T_1\ge 4.  \lb{5.53}\ee

On the other hand, because $i(c^{T_1})=i(c^{T})=R$, replacing $T$ by $T_1$ equalities in
(\ref{5.52}) still hold, and then it yields
$$  T_1 - 3(\{T_1\sg_1\}+\{T_1\sg_2\}) = T - 3(\{T\sg_1\}+\{T\sg_2\}) = T-3.  $$
Together with (\ref{5.53}), it implies that
\be  0 < 3(\{T_1\sg_1\}+\{T_1\sg_2\}) = 3 + (T_1-T) \le 3-4 = -1. \lb{5.54}\ee
This contradiction proves $k_2=0$. Then Claim 3 for $m\in\N$ follows from Lemma 2.2.

\medskip

Because all $i(c^m)$ are odd, by (\ref{5.42}), (\ref{5.43}), Claim 3 and
Lemma 2.1 we obtain that $c^m$s with $m\ge T+1$ have no contributions to
$M_j$s with $0\le j\le R+1$, and $c^m$s with $1\le m\le T$ have only
contributions to $M_j$s with $0\le j\le R+1$. More precisely,
$\sum_{2j-1=1}^RM_{2j-1}$ is completely contributed by $c^m$s with
all integer $m\in [1,T]$ which are not in $n\N$ and each $c^m$ contributes a $1$,
as well as by $c^{mn}$s with all integer $m\in [1,T/n]$ and each $c^{mn}$ contributes
a $k_0$. And $\sum_{2j=0}^{R+1}M_{2j}$ is completely contributed by $c^{mn}$s with
$m\in\N$ satisfying $mn\in [n,T]$ and each $c^{mn}$ contributes a $k_1$. Thus we
have
\bea  \sum_{j=0}^{R+1}(-1)^j M_j
&=& \sum_{2j=0}^{R+1}M_{2j}-\sum_{2j-1=1}^RM_{2j-1}   \nn\\
&=& \frac{T}{n}k_1-\left(\frac{T}{n}k_0+T-\frac{T}{n}\right)  \nn\\
&=& \frac{k_1-(k_0+n-1)}{n}T.    \lb{5.55}\eea
On the other hand, by Claim 3 and Lemma 2.4, we have
\be  \frac{k_1-(k_0+n-1)}{n\hat{i}(c)} = -\frac{3}{2}.  \lb{5.56}\ee
Thus by (\ref{5.55}), (\ref{5.56}), the Morse inequality and (\ref{4.37}), we obtain
\be -\frac{3T}{2}\hat{i}(c) = \sum_{j=0}^{R+1}(-1)^jM_j
       \ge \sum_{j=0}^{R+1}(-1)^jb_j = -\sum_{2k-1=1}^R b_{2k-1}=-\frac{3(R-1)}{2}.
           \lb{5.57}\ee
It implies
\be   R-1 \ge T\hat{i}(c).  \lb{5.58}\ee
But on the other hand, by (\ref{5.29}) we have
\bea R-T\hat{i}(c)
&=& (-2T+2\sum_{j=1}^3E(T\sg_j)-3) - T(-2+2\sum_{j=1}^3\sg_j)   \nn\\
&=& 2\sum_{j=1}^2(E(T\sg_j)-T\sg_j) - 3     \nn\\
&=& 2\sum_{j=1}^2 (1-\{T\sg_j\}) - 3    \nn\\
&<& 1,   \lb{5.59}\eea
which contradicts to (\ref{5.58}) and completes the proofs in this subcase and Case 2-1.

\medskip

{\bf Case 2-2:} {\it $G= N_1(-1,1)$.}

In this case, by $i(c)=1$ and Theorem 3.3, we have the iteration formula
\be  i(c^m) = - m + 2\sum_{j=1}^2E(m\sg_j) - 2, \qquad
        \nu(c^m) = \frac{1+(-1)^m}{2},  \qquad  \forall\,m\ge 1.  \lb{5.60}\ee
Then we have $n=n(c)=2$ and
\be i(c^m)=m \quad(\mod\;2), \qquad\nu(c^{2m-1})=0\;\;\mbox{and}\;\;\nu(c^{2m})=1,
     \qquad \forall\; m\in\N. \lb{5.61}\ee

By Lemma 2.2, it yields
$k_0^-(c^{2k})=k_0^-(c^2)=0$ for all $k\in\N$. Because the
iterates $c^{2k-1}$ with $k\in\N$ contribute only to the odd-th
Morse-type numbers, we obtain
\be M_{2k}=\,^{\#}\{m\in\N\,|\,i(c^m)=2k\}k_0^-(c^2)=0. \lb{5.62}\ee
Together with the Morse inequality, it implies that for any $k\ge 1$,
\be  b_{2k-1} = M_{2k-1}\nn\\
  = \,^{\#}\{j\in\N\,|\, i(c^{2j})=2k-2\}k_1^-(c^2)
            + \,^{\#}\{j\in\N\,|\, i(c^{2j-1})=2k-1\}. \lb{5.63}\ee
Note that $b_1=1$ when $d=2$, and $b_1=0$ when $d=4$. By the facts $i(c)=1$,
$\nu(c)=0$ and (\ref{5.63}) in this case we must have
\be   d=h=2.   \lb{5.64}\ee
Therefore (\ref{5.15}) holds again by Lemma 2.6.

Let $k_1\equiv k_1^-(c^2)\in\{0,1\}$. Then by Lemma 2.4, we obtain
$$ \frac{-1-k_1}{2(2(\sg_1+\sg_2)-1)}
  = \sum_{1\le m\le 2,\,0\le l\le 1}\frac{(-1)^{i(c^m)+l}k_l^{\ep}(c^m)}{n\hat{i}(c)}
  = B(2,2) = - \frac{3}{2},   $$
which yields
\be   \sg_1 + \sg_2 = \frac{4+k_1}{6}.  \lb{5.65}\ee

{\bf Claim 4.} {\it $k_1^-(c^{2m})=k_1=1$ for all $m\in\N$.}

In fact, assume $k_1=0$. Then all $i(c^{2k})$ with $k\ge 1$ have no contribution
to the odd-th Morse-type number $M_{2j-1}$ with $j\in\N$. In addition, by (\ref{5.65}),
we obtain $3(\sg_1+\sg_2) = 2$. Because both of $3\sg_1$ and $3\sg_2$ are irrational,
it yields $[3\sg_1]+[3\sg_2]=1$. Thus by (\ref{5.60}), we obtain $i(c^3)=1$. Together
with $i(c)=1$, it yields $M_1\ge 2$. It contradicts to the fact $M_1=b_1=1$ obtained
from (\ref{5.62}), (\ref{5.63}), (\ref{5.64}) and Lemma 2.6. By Lemma 2.2, Claim 4
is proved.

Next we estimate $i(c^m)$ using Lemma 5.1. By Claim 4, (\ref{5.65}) becomes
\be   \sg_1 + \sg_2 = \frac{5}{6}.  \lb{5.66}\ee
Then by Lemma 5.1 we obtain
\be  [\frac{5m}{6}]-1 \le [m\sg_1] + [m\sg_2] \le [\frac{5m}{6}], \qquad \forall\; m\in\N. \lb{5.67}\ee

Thus by (\ref{5.60}) and (\ref{5.9}) for $m=6k\in\N$ we obtain
\be  i(c^{6k}) = -6k + 2([6k\frac{5}{6}]-1) + 2 = 4k,  \quad \forall\;k\in\N.  \lb{5.68}\ee

By (\ref{5.60}) and (\ref{5.67}) for $m=6k+1\in\N$, we obtain
$$ -(6k+1)+2([(6k+1)\frac{5}{6}]-1)+2 \le i(c^{6k+1}) \le -(6k+1)+2[(6k+1)\frac{5}{6}]+2. $$
That is,
\be  4k-1 \le i(c^{6k+1}) \le 4k+1, \qquad \forall\;k\in\N_0.  \lb{5.69}\ee

Similarly for $m=6k+2\in\N$, we obtain
$$ -(6k+2)+2([(6k+2)\frac{5}{6}]-1)+2 \le i(c^{6k+2}) \le -(6k+2)+2[(6k+2)\frac{5}{6}]+2. $$
It yields
\be  4k \le i(c^{6k+2}) \le 4k+2, \qquad \forall\;k\in\N_0.  \lb{5.70}\ee

For $m=6k+3\in\N$, we obtain
$$ -(6k+3)+2([(6k+3)\frac{5}{6}]-1)+2 \le i(c^{6k+3}) \le -(6k+3)+2[(6k+3)\frac{5}{6}]+2. $$
It yields
\be  4k+1 \le i(c^{6k+3}) \le 4k+3, \qquad \forall\;k\in\N_0.  \lb{5.71}\ee

For $m=6k+4\in\N$, we obtain
$$ -(6k+4)+2([(6k+4)\frac{5}{6}]-1)+2 \le i(c^{6k+4}) \le -(6k+4)+2[(6k+4)\frac{5}{6}]+2. $$
It yields
\be  4k+2 \le i(c^{6k+4}) \le 4k+4, \qquad \forall\;k\in\N_0.  \lb{5.72}\ee

For $m=6k+5\in\N$, we obtain
$$ -(6k+5)+2([(6k+5)\frac{5}{6}]-1)+2 \le i(c^{6k+5}) \le -(6k+5)+2[(6k+5)\frac{5}{6}]+2. $$
It yields also
\be  4k+3 \le i(c^{6k+5}) \le 4k+5, \qquad \forall\;k\in\N_0.  \lb{5.73}\ee

Then using similar arguments in the proof of Claim 1, we have

\medskip

{\bf Claim 5:} {\it $i(c)=1$, $i(c^2)=2$, $i(c^3)=3$, $i(c^4)=4$, $i(c^5)=5$,
$i(c^6)=4$, $i(c^7)\le 5$.}

\medskip

In fact, (\ref{5.15}) is crucial in the following. Note that $c$ contributes a $1$
to $M_1=b_1=1$ by the facts $i(c)=1$ and $\nu(c)=0$. Thus $i(c^m)\ge 2$ for all $m\ge 2$.
Then by Claim 4 and (\ref{5.70})-(\ref{5.71}) with $k=0$ we obtain $i(c^2)\le 2$ and
$i(c^3)\le 3$. Thus by (\ref{5.61}) we obtain $i(c^2)=2$ and $i(c^3)=3$.

By (\ref{5.68}) with $k=1$ and (\ref{5.72})-(\ref{5.73}) with $k=0$, we obtain
$i(c^6)=4$, $i(c^4)\le 4$, $i(c^5)\le 5$. By Claim 4 and (\ref{5.15}) we obtain
$i(c^4)=4$ and $i(c^5)=5$.

Then by (\ref{5.69}) with $k=1$ we obtain $i(c^7)\le 5$. Claim 5 is proved.

\medskip

Now by Claims 4 and 5, each $i(c^k)$ with $1\le k \le 7$ contributes $1$ to the
Morse-type numbers $M_1+M_3+M_5$. Thus by (\ref{5.15}), (\ref{5.62})-(\ref{5.63})
we obtain
\be  6 = \sum_{j=0}^5b_j = \sum_{j=0}^5M_j\ge 7.  \lb{5.74}\ee
Contradiction!

\medskip

{\bf Case 2-3:} {\it $G=N_1(1,a)$ with $a=0$ or $1$.}

\medskip

In this case, by $i(c)=1$ and Theorem 3.3 we have the formula
\be i(c^m)=2\sum_{j=1}^2[m\sg_j]+1, \qquad \nu(c^m)=2-a, \qquad\forall\,m\ge 1.  \lb{5.75}\ee
Note that all $i(c^m)$ with $m\ge 1$ are odd and non-decreasing in $m$. Because $b_1=1$ when
$d=h=2$, and $b_3=1$ when $d=4$ and $h=1$, to generate the non-zero Morse-type number
$M_1\ge b_1$ or $M_3\ge b_3$, there must hold $k_0(c)+k_2^+(c)=1$ and $k_1^+(c)=0$. Thus by
the Morse inequality and Lemmas 2.2 and 2.4, it yields
\bea
&& M_{2k}=0,\quad M_{2k-1}=b_{2k-1},\,\forall\,k\in\N_0,   \lb{5.76}\\
&& \sg_1 + \sg_2 = \frac{-1}{2B(d,h)}.  \lb{5.77}\eea

If $d=4$ and $h=1$, we have $B(d,h)=-2/3$ and $b_3=b_5=b_7=1<2=b_9$ by Lemma 2.5.
In order to get $M_3=b_3=1$, by $i(c)=1$ and Lemma 2.2 it yields $k_2^+(c^m)=k_2^+(c)=1$
and $k_0(c^m)=k_0(c)=k_1^+(c^m)=k_1^+(c)=0$ for all $m\ge 1$. Note that by (\ref{5.77}) we have
$$   \sg_1+\sg_2 = \frac{3}{4}.  $$
Thus $[4\sg_1]+[4\sg_2]=3-1=2$ by Lemma 5.1. So we have $i(c^4)=5$ by (\ref{5.75}). Thus
we obtain
$$  4\le \sum_{2j-1=3}^7M_{2j-1}=\sum_{2j-1=3}^7b_{2j-1}=3,  $$
a contradiction.

If $d=h=2$, we have $B(d,h)=-3/2$ and $b_1=1$ by Lemma 2.5. To generate $M_1=b_1=1$, we should
have $1=k_0(c)=k_0(c^m)$ and $0=k_1^+(c)=k_2^+(c)=k_1^+(c^m)=k_2^+(c^m)$ for all $m\in\N$ by
$i(c)=1$ and Lemma 2.2. Notice that
$$   \sg_1+\sg_2=\frac{1}{3}  $$
holds by (\ref{5.77}). Then $[3\sg_1]+[3\sg_2]=1-1=0$ by Lemma 5.1. Thus $i(c^3)=1$ by
(\ref{5.75}). Then by the monotone increasing of the Morse indices $i(c^m)$ in $m$, we
obtain $3\le M_1=b_1=1$, a contradiction.

\medskip

This completes the proof of Step 2 for $i(c)=1$.

\medskip

{\bf Step 3.} {\it $i(c)\ge 2$.}

\medskip

Because $i(c^m)\ge i(c)$ for $m\ge 1$ due to the Bott formula, to generate the non-trivial
homology $H_{d-1}(\ol{\Lm}M,\ol{\Lm}^0 M;\Q)$, then the Morse index of $c$ must satisfy
$i(c)\le d-1$. Thus we must have $d=4$ and $h=1$. In other words, the manifold is rationally
homotopic to $S^4$. We continue the proof in two cases according to the value of $i(c)$.

\medskip

{\bf Case 3-1:} {\it $i(c)=2$.}

\medskip

By Proposition 3.4, we must have $G=N_1(1,-1)$ in (\ref{5.4}).
Thus, by Theorem 3.3, we have
\be  i(c^m) = 2([m\sg_1]+[m\sg_2]) + 2 \quad\mbox{and}\quad \nu(c^m)=1, \qquad
       \forall\;m\in\N. \lb{5.78}\ee
Thus in this case, we have $i(c^m)\in 2\Z$ for all $m\in\N$ and are non-decreasing in $m$,
and then $n=n(c)=1$. Thus by Lemma 2.4 we have the identity
\be  -\frac{3}{2}(k_0(c)-k_1^+(c)) = 2(\sg_1+\sg_2) = \hat{i}(c)>0, \lb{5.79}\ee
which implies $k_1^+(c^m)=k_1^+(c)=1$ and $k_0(c^m)=0$ for all $m\ge 1$ by Lemma 2.2.
So (\ref{5.79}) becomes
\be   \sg_1 + \sg_2 = \frac{3}{4}.  \lb{5.80}\ee
By (\ref{5.9}) in Lemma 5.1 we obtain $[4\sg_1]+[4\sg_2]=3-1=2$, and then
\be  i(c^4)= 2([4\sg_1]+[4\sg_2]) + 2 = 6.  \lb{5.81}\ee
Thus by Theorem 3.13 we get
\be  i(c^m)\le 6, \qquad \forall\;m= 1, 2, 3, 4.  \lb{5.82}\ee

From the above discussion, for all integer $k\ge 0$ we get
\bea
M_{2k} &=& \;^\#\{m\ge1:i(c^m)=2k)\}\,k_0(c) = 0,  \lb{5.83}\\
M_{2k+1} &=& \;^\#\{m\ge1:i(c^m)=2k\}.  \lb{5.84}\eea
Thus we have $M_1+M_3+M_5+M_7\ge 4$ by (\ref{5.82}) and (\ref{5.84}). Thus the Morse
inequality and Lemma 2.5 again yield a contradiction:
\be  -4\ge \sum_{q=0}^8 (-1)^q M_q \ge \sum_{q=0}^8 (-1)^q b_q = -3.  \lb{5.85}\ee

\medskip

{\bf Case 3-2:} {\it $i(c)\ge 3$.}

\medskip

Note that by Theorem 3.13 it yields $i(c^{m+1})\ge i(c^m)$ for all $m\ge 1$.

By Lemma 2.4, both $\sg_1$ and $\sg_2$ are linearly dependent over $\Q$. Thus
we must have $A=1$ in Theorem 3.21 and there exists some $T\in 12n\N$ with
$n=n(c)$ being the analytical period of $c$ such that
\bea
i(c^m)-i(c^T) &\ge& i(c)+p_0+p_-+ (q_0+q_+) + r-2\equiv \xi(c), \qquad \forall\;m\ge T+1,  \lb{5.86}\\
i(c^T)-i(c^m) &\ge& i(c)-r + p_-+p_0 + k-(q_0+q_+) \ge 0, \qquad \forall\; 1\le m\le T-1, \lb{5.87}\eea
where we used the fact $k=2$ in Theorem 3.21.

Let $\tau^-(m)=\frac{1-(-1)^m}{2}$ for any $m\in\N$. Note that
\be \nu(c^n)=2(p_0+q_0)+2(r-2)+p_-+q_++p_++q_-. \lb{5.88}\ee
From $i(c)\ge 3$ and the fact $r-2+p_0+q_0+p_++q_-\le 1$, we get
$$ \xi(c) = i(c) + \nu(c^n) - (r-2+p_0+q_0+p_++q_-) \ge \nu(c^n) + 1 + \tau^-(i(c^T)+\nu(c^n)). $$
Then (\ref{5.86}) becomes
\be i(c^m)-i(c^T) \ge \nu(c^n) + 1 + \tau^-(i(c^T)+\nu(c^n)), \qquad \forall\;m\ge T+1.  \lb{5.89}\ee

Let $R=i(c^T)$, $\bar{\nu}=\nu(c^T)=\nu(c^n)$ and
$\tdR \equiv R + \bar{\nu} + \tau^-(R+\bar{\nu})\in 2\Z$.
It follows from (\ref{5.87}) that all iterations $c^m$ with $1\le m\le T$ contribute only to
the Morse-type numbers $M_q$ for $0\le q\le R+\bar{\nu}$, and from (\ref{5.89}) that
all the iterations $c^m$ with $m\ge T+1$ do not contribute to these Morse-type numbers
$M_q$ with $0\le q\le R+\bar{\nu}$. Thus it yields
\bea  \sum_{q=0}^{R+\bar{\nu}}(-1)^q M_q
&=& \sum_{0\le q\le R+\bar{\nu}\atop 1\le m\le T}(-1)^q\dim\ol{C}_q(E,c^m)  \nn\\
&=& \sum_{m=1}^T\left(\sum_{q=0}^{R+\bar{\nu}}(-1)^{i(c^m)+(q-i(c^m))}k_{q-i(c^m)}^{\ep(c^m)}(c^m)\right)\nn\\
&=& \sum_{m=1}^T\left(\sum_{q=0}^{i(c^m)+\nu(c^m)}
           (-1)^{i(c^m)+(q-i(c^m))}k_{q-i(c^m)}^{\ep(c^m)}(c^m)\right)\nn\\
&=& \sum_{m=1}^T\left(\sum_{l_m=0}^{\nu(c^m)}(-1)^{i(c^m)+l_m}k_{l_m}^{\ep(c^m)}(c^m)\right)\nn\\
&=& \frac{T}{n}\sum_{1\le m\le n \atop 0\le l_m\le\nu(c^m)}(-1)^{i(c^m)+l_m}k_{l_m}^{\ep(c^m)}(c^m)\nn\\
&=& T\hat{i}(c)B(4,1),  \lb{5.90}\eea
where we used (\ref{5.87}) and (\ref{5.89}) in the first equality, Lemma 2.1 in the second one, (i) of
Lemma 2.2 in the third and fourth ones, Lemma 2.3 in the fifth one, and Lemma 2.4 in the sixth one.

In this case, by Lemma 2.5 only $b_q$s with odd $q\ge 3$ are non-zero. By (\ref{5.90}),
the Morse inequality and Lemma 2.5 we obtain
\be  T\hat{i}(c)B(4,1) = \sum_{j=0}^{\tdR}(-1)^j M_j
  \ge \sum_{j=0}^{\tdR} (-1)^jb_{j} = -\sum_{2q-1=1}^{\tdR-1} b_{2q-1}
  \ge \frac{5-2\tdR}{3}.   \lb{5.91}\ee
Here $M_{R+\bar{\nu}+\tau^-(R+\bar{\nu})}=0$ by (\ref{5.87}) and (\ref{5.89}) when $R+\bar{\nu}$
is odd. This fact is used in the first equality in (\ref{5.91}) when $R+\bar{\nu}$ is odd.
Note also that in the first inequality of (\ref{5.91}), the evenness of $\tdR$ implies the
availability of the Morse inequality.

On the other hand, by Lemma 2.4 we obtain
$$  \frac{-t}{n(s+2(\sg_1+\sg_2)+\frac{q}{p})} = B(4,1) = -\frac{2}{3},  $$
for some integers $s$, $t$, $q$ and $p$, where we write $q/p=\th_3/\pi\in [0,1)\cap\Q$ with
$(p,q)=1$ when $q>0$ for the possible term $R(\th_3)$. From this identity we obtain
$$  \sg_1 + \sg_2 = \frac{3t}{4n} - \frac{s}{2} - \frac{q}{2p} = \frac{b}{4n}  $$
for some integer $b>0$, where we have used the fact $p|n$ which follows from Definition 3.6
of $n=n(c)$. Thus according to the choice of $T\in 12n\N$, by (\ref{5.10}) we obtain
\be  \{T\sg_1\}+\{T\sg_2\}=1.  \lb{5.92}\ee
Also note that (\ref{5.88}) yields
\be  \nu(c^n)-(r-2+p_-+p_0+q_++q_0)=p_0+q_0+p_++q_-+r-2\le 1. \lb{5.93}\ee

By Theorem 3.3 and (\ref{5.92})-(\ref{5.93}) we obtain that
the integer on the right hand side of (\ref{5.90}) satisfies
\bea T\hat{i}(c)B(4,1)
&=& -\frac{2T}{3}\left(i(c)+p_-+p_0-r+2\sum_{j=1}^r\sg_j\right) \nn\\
&=& -\frac{2}{3}\left(T(i(c)+p_-+p_0-r)
      + 2\sum_{j=1}^rE\left(T\sg_j\right) -2\right)  \nn\\
&=& -\frac{2}{3}\left(i(c^T)+r+p_-+p_0+q_++q_0-2\right)  \nn\\
&\le& -\frac{2}{3}(R+\nu(c^n)-1)  \nn\\
&\le&  -\frac{2}{3}(\tdR-2),   \lb{5.94}\eea
where the first equality follows from (\ref{3.9}), the second equality follows from
(\ref{5.92}) and the fact
$$  \sum_{j=1}^rT\sg_j = \sum_{j=1}^r([T\sg_j]+\{T\sg_j\})
      = \sum_{j=1}^rE(T\sg_j)-2 + \sum_{j=1}^2\{T\sg_j\} = \sum_{j=1}^rE(T\sg_j)-1,  $$
the third equality follows from (\ref{3.7}) with $m=T\in 2\N$, the first inequality
follows from (\ref{5.93}), and the last inequality follows from the definition of
$\tdR$.

Now (\ref{5.91}) and (\ref{5.94}) yield a contradiction.

\medskip

The proof of Theorem 1.2 is complete. \hfill\hb

\medskip

\setcounter{equation}{0}
\section{On compact simply connected reversible Finsler manifolds}

In this section, we study closed geodesics on compact simply connected reversible
Finsler manifolds, including Riemannian manifolds, and give the proofs of the main
Theorems 1.1 and 1.3 about closed geodesics on $4$-dimensional compact simply connected
reversible Finsler manifolds.

For any reversible Finsler as well as Riemannian metric $F$ on a compact
manifold $M$, the energy functional $E$ is symmetric on every loop
$f\in \Lm M$ and its inverse curve $f^{-1}$ defined by $f^{-1}(t)=f(1-t)$. Thus
these two curves have the same energy $E(f)=E(f^{-1})$ and play the same roles
in the variational structure of the energy functional $E$ on $\Lm M$.
Specially, the $m$-th iterates $c^m$ and $c^{-m}$ of a closed
geodesic $c$ and its inverse curve $c^{-1}$ have precisely the same
Morse indices, nullities, and critical modules. Let $n=n(c)$. So there holds
\be   \dim\ol{C}_*(E,c^m)=\dim\ol{C}_*(E,c^{-m}).    \lb{6.1}\ee
Thus if $c$ is the only geometrically distinct prime closed geodesic on $M$,
then all the Morse type numbers must be even, i.e.,
\be   M_j \in 2\N_0, \qquad \forall\; j\in\Z,  \lb{6.2}\ee
and the identity in Lemma 2.4 becomes
\be  2\sum_{0\le l_m\le \nu(c^m)\atop 1\le m\le n}
  (-1)^{i(c^m)+l_m}k_{l_m}^{\ep(c^m)}(c^m)=n\hat{i}(c)B(d,h). \lb{6.3}\ee
From this consideration we get the following result.

\medskip

{\bf Theorem 6.1.} {\it Theorems 4.3 and 4.4 hold for reversible Finsler (as well as
Riemannian) metric on the corresponding manifold $(M,F)$ too. Therefore Theorem 1.1
holds. }

{\bf Proof.} The current version of Theorem 4.3 works in the reversible Finsler
metric case without any changes by the same reason as we have explained in Remark 7.1
of \cite{LoD1}. Note that now the integer $\ka$ in (\ref{4.8}) is even by the above
reason.

For the Claim 1 of Theorem 4.4 with a reversible Finsler metric on $M$, by the same
reason, the above proof of Theorem 4.4 works without any change and shows that the
only geometrically distinct prime closed geodesic $c$ which can not be rational in
the reversible case.

For Theorem 4.4 with a reversible Finsler metric on $M$ and only one geometrically
distinct prime closed geodesic $c$ on $M$ which is completely non-degenerate, the
above proof of Claim 2 in Theorem 4.4 with minor modifications works too. In fact,
Lemma 4.1 and (\ref{6.2}) yields a much simpler proof, because we get the following
contradiction immediately
\be  1 = b_{dh-1} = M_{dh-1} \in 2\N_0,  \lb{6.4}\ee
where $dh=\dim M$. Therefore Claim 2 of Theorem 4.4 holds too in the reversible
Finsler metric case.

Thus Theorem 1.1 holds. \hfill\hb

\medskip

Now we can give

{\bf The proof of Theorem 1.3.} This proof is similar to that of Theorem 1.2 in
Section 5. Next we follow the classification used in the proof of Theorem 1.2
and indicate only some necessary changes and omit the details.

{\bf Step 1.} $i(c)=0$.

Following the study in Step 1 of the proof of Theorem 1.2, we have $i(c)=0$ and
$G=N_1(-1,1)$ in (\ref{5.4}). By (\ref{6.1}), the positive numbers $1/|B(d,h)|$
should be replaced by $2/|B(d,h)|$ in (\ref{5.12}). Then similar arguments yield
(\ref{5.15}), specially by (\ref{6.2}) we obtain the following contradiction
\be  1 = b_1 = M_1 \in 2\N_0,  \lb{6.5}\ee
and then complete the proof in Step 1.

{\bf Step 2.} $i(c)=1$.

{\sl Case 2-1.} As in the proof of Theorem 1.2, we distinguish two subcases.

{\sl Subcase 2-1-1.} {\it $d=4$ and $h=1$. }

Replacing $k_1^+(c^n)$ by $2k_1^+(c^n)$ in (\ref{5.33}), by the same proof we
get Claim 2. Then replacing $c^m$ by $c^m$ and $c^{-m}$, $c^{mn}$ by $c^{mn}$
and $c^{-mn}$ in the paragraph below (\ref{5.36}), instead of (\ref{5.37}) and
(\ref{5.38}), by (\ref{6.2}) and (\ref{6.3}) we obtain
\bea
&& \sum_{j=0}^{R+1}(-1)^j M_j = 2T\frac{k_1-n+1}{n},   \lb{6.6}\\
&& 2\frac{k_1-(n-1)}{n\hat{i}(c)} = -\frac{2}{3}.  \lb{6.7}\eea
Then using (\ref{6.6}) and (\ref{6.7}), the same proofs from (\ref{5.39}) to
(\ref{5.41}) yield a contradiction.

{\sl Subcase 2-1-2.} $d=h=2$.

In this subcase, note that we have still (\ref{5.46}) if $k_2^+(c^n)=1$.
Thus the contradiction $1=b_1=M_1\in 2\N_0$ yields Claim 3.

Now as in the above Subcase 2-1-1, the proofs in (\ref{5.55}) to (\ref{5.59})
yield a contradiction.

{\sl Case 2-2.} Similarly to (\ref{5.60})-(\ref{5.64}), we obtain
$1=b_1=M_1\in 2\N_0$, contradiction!

{\sl Case 2-3.} In this case, from (\ref{5.76})-(\ref{5.77}) and (\ref{6.2}) we
obtain the contradiction $1=b_{d-1}=M_{d-1}\in 2\Z$ with $d=2$ or $d=4$.

{\bf Step 3.} $i(c)\ge 2$.

{\sl Case 3-1.} $i(c)=2$.

In this case $G=N_1(1,-1)$. By (\ref{6.3}) the identity
(\ref{5.80}) now becomes
\be  -3(k_0(c)-k_1^+(c)) = 2(\sg_1+\sg_2) = \hat{i}(c)>0,  \lb{6.8}\ee
with $k_0(c)=0$ and $k_1^+(c)=1$, and
\be   \sg_1 + \sg_2 = \frac{3}{2}.  \lb{6.9}\ee
By (\ref{5.79}), this specially implies
$$  M_{2k}=0, \qquad M_{2k+1}=b_{2k+1} \qquad \forall\; k\in \N_0. $$
Thus by Lemma 5.1 we obtain $[2\sg_1]+[2\sg_2]=3-1=2$, and then
\be  i(c^2) = 2([2\sg_1]+[2\sg_2]) + 2 = 6.  \lb{6.10}\ee
Thus by the monotone increasing of $i(c^m)$ in $m$ from Theorem 3.13, we
obtain the following contradiction
$$  0 = M_5 = b_5 \ge 1.  $$

{\sl Case 3-2.} $i(c)\ge 3$.

Because the contributions of $c^{-m}$ with $m\ge 1$, similarly to (\ref{5.90})
by (\ref{6.2}) and (\ref{6.3}) we obtain
\bea  \sum_{q=0}^{R+\bar{\nu}}(-1)^q M_q
&=& 2\sum_{0\le q\le R+\bar{\nu}\atop 1\le m\le T}(-1)^q\dim\ol{C}_q(E,c^m)  \nn\\
&=& \frac{2T}{n}\sum_{1\le m\le n \atop 0\le l_m\le\nu(c^m)}(-1)^{i(c^m)+l_m}k_{l_m}^{\ep(c^m)}(c^m)\nn\\
&=& T\hat{i}(c)B(4,1).  \lb{6.11}\eea
Then by the same proof of (\ref{5.91}) and (\ref{5.94}) we obtain a contradiction.

The proof of Theorem 1.3 is complete. \hfill\hb

\medskip

{\bf Acknowledgements.} The authors thank sincerely the referee for his/her
careful reading of the manuscript and valuable comments.

\bibliographystyle{abbrv}

\end{document}